\RequirePackage{amsthm}%
 
\documentclass[sn-apa]{sn-jnl}

\usepackage{graphicx}%
\usepackage{amsmath,amssymb,amsfonts}%
\usepackage{amsthm}%
\usepackage{mathrsfs}%
\usepackage[title]{appendix}%
\usepackage{xcolor}%
\usepackage{textcomp}%
\usepackage{manyfoot}%
\usepackage{booktabs}%
\usepackage{algorithm}%
\usepackage{algorithmicx}%
\usepackage{algpseudocode}%
\usepackage{listings}%

\usepackage{multirow}
\usepackage{makecell}
\usepackage{bm}
\usepackage{enumitem}
\usepackage{adjustbox}
\usepackage{anyfontsize}
\usepackage{soul}

\usepackage{subfigure}

\newcolumntype{L}[1]{>{\raggedright\let\newline\\\arraybackslash\hspace{0pt}}m{#1}}
\newcolumntype{C}[1]{>{\centering\let\newline\\\arraybackslash\hspace{0pt}}m{#1}}
\newcolumntype{R}[1]{>{\raggedleft\let\newline\\\arraybackslash\hspace{0pt}}m{#1}}

\theoremstyle{thmstyleone}%
\newtheorem{theorem}{Theorem}

\theoremstyle{thmstyletwo}%

\theoremstyle{thmstylethree}%

\def\cA{\mathcal A}

\def\cD{\mathcal D}

\def\cF{\mathcal F}
\def\cG{\mathcal G}
\def\cH{\mathcal H}

\def\cN{\mathcal N}

\def\cQ{\mathcal Q}

\def\cZ{\mathcal Z}

\newcommand{\bd}{{\bf d}}

\newcommand{\bff}{{\bf f}}
\newcommand{\bg}{{\bf g}}

\newcommand{\bh}{{\bf h}}

\newcommand{\bj}{{\bf j}}

\newcommand{\bp}{{\bf p}}

\newcommand{\bt}{{\bf t}}

\newcommand{\bv}{{\bf v}}

\newcommand{\bx}{{\bf x}}
\newcommand{\bX}{{\bf X}}
\newcommand{\by}{{\bf y}}
\newcommand{\bY}{{\bf Y}}
\newcommand{\bz}{{\bf z}}
\newcommand{\bZ}{{\bf Z}}

\newcommand{\bbI}{{\mathbb I}}
\newcommand{\bbN}{{\mathbb N}}
\newcommand{\bbP}{{\mathbb P}}

\newcommand{\bbR}{{\mathbb R}}

\newcommand{\E}{\mathbb{E}}

\newcommand{\bbeta}{\bm{\beta}}

\newcommand{\bepsilon}{\bm{\epsilon}}

\newcommand{\ie}{{\it i.e.}}
\newcommand{\cf}{{\it cf.}}
\newcommand{\eg}{{\it e.g.}}

\newcommand{\iid}{{i.i.d.}}

\newcommand{\Holder}{H\"{o}lder }

\newcommand{\bc}{\begin{center}}
\newcommand{\ec}{\end{center}}
\newcommand{\be}{\begin{equation}}
\newcommand{\ee}{\end{equation}}
\newcommand{\ba}{\begin{array}}
\newcommand{\ea}{\end{array}}
\newcommand{\bean}{\setlength\arraycolsep{1pt}\begin{eqnarray*}}
\newcommand{\eean}{\end{eqnarray*}}
\newcommand{\bea}{\setlength\arraycolsep{1pt}\begin{eqnarray}}
\newcommand{\eea}{\end{eqnarray}}
\newcommand{\ben}{\begin{enumerate}}
\newcommand{\een}{\end{enumerate}}
\newcommand{\bed}{\begin{itemize}}
\newcommand{\eed}{\end{itemize}}

\DeclareMathOperator*{\argmax}{argmax}

\DeclareMathOperator*{\minimize}{minimize}

\newcommand{\bzero}{{\bf 0}}

\newcommand{\Id}{\bbI}
\newcommand{\vertiii}[1]{{\left\vert\kern-0.25ex\left\vert\kern-0.25ex\left\vert #1 
    \right\vert\kern-0.25ex\right\vert\kern-0.25ex\right\vert}}

\begin{document}

\title[Convergence rates for GAN]{Rates of convergence for nonparametric estimation of singular distributions using generative adversarial networks}


\author[1]{\fnm{Jeyong} \sur{Lee}}

\author[1]{\fnm{Hyeok Kyu} \sur{Kwon}}

\author*[1]{\fnm{Minwoo} \sur{Chae}}\email{mchae@postech.ac.kr}

\affil[1]{\orgdiv{Department of Industrial and Management Engineering}, \orgname{Pohang University of Science and Technology}, \orgaddress{\city{Pohang}, \country{South Korea}}}

\abstract{It is common in nonparametric estimation problems to impose a certain low-dimensional structure on the unknown parameter to avoid the curse of dimensionality. This paper considers a nonparametric distribution estimation problem with a structural assumption under which the target distribution is allowed to be singular with respect to the Lebesgue measure. In particular, we investigate the use of generative adversarial networks (GANs) for estimating the unknown distribution and obtain a convergence rate with respect to the $L^1$-Wasserstein metric. The convergence rate depends only on the underlying structure and noise level. More interestingly, under the same structural assumption, the convergence rate of GAN is strictly faster than the known rate of VAE in the literature. We also obtain a lower bound for the minimax optimal rate, which is conjectured to be sharp at least in some special cases. Although our upper and lower bounds for the minimax optimal rate do not match, the difference is not significant.}


\maketitle

\addtocontents{toc}{\protect\setcounter{tocdepth}{-1}}

\section{Introduction}

Given $D$-dimensional observations $\bX_1, \ldots, \bX_n$ following $P_0$, suppose that we are interested in inferring the underlying distribution $P_0$ or related quantities such as its density function or the manifold on which $P_0$ is supported.
The inference of $P_0$ is fundamental in unsupervised learning, and there are numerous inferential methods available in the literature. We refer to Chapter 14 of \cite{hastie2009elements} for various methods.

In this paper, $\bX_i$ is modeled as $\bX_i = \bg(\bZ_i) + \bepsilon_i$ for some function $\bg: \cZ \to \bbR^D$.
Here, $\bZ_i$ is a latent variable following the known distribution $P_Z$ supported on $\cZ \subset \bbR^d$, and $\bepsilon_i$ is an error following a normal distribution $\cN(\bzero_D, \sigma^2 \Id_D)$, where $\bzero_D$ and $\Id_D$ denote a $D$-dimensional vector of zeros and an identity matrix, respectively.
The dimension $d$ of the latent variable $\bZ_i$ is typically much smaller than $D$.
This model is often called a (non-linear) \emph{factor model} in statistical communities (\cite{yalcin2001nonlinear, kundu2014latent}) and a \emph{generative model} in machine learning societies \citep{goodfellow2014generative, kingma2013auto}.
Throughout the paper, we use the latter terminology.
Accordingly, $\bg$ will be referred to as a \emph{generator}.

Recent advances in deep learning have expanded the use of generative models by modeling $\bg$ through deep neural networks (DNN), also known as \emph{deep generative models}.
Two approaches are popularly used for estimating $\bg$.
Variational autoencoder (VAE; \cite{kingma2013auto,rezende2014stochastic}) is perhaps the most well-known algorithm for constructing an estimator $\hat\bg$ using a likelihood approach.
The other approach is \emph{generative adversarial networks (GAN)}. Originally proposed by \cite{goodfellow2014generative}, GAN has been extended in several directions. One of its extensions considers general integral probability metrics (IPM) as loss functions.
Sobolev GAN \citep{mroueh2017sobolev}, maximum mean discrepancy GAN \citep{li2017mmd} and Wasserstein GAN \citep{arjovsky2017wasserstein} are important examples in this direction.
Another important direction of generalization is the development of novel architectures for generators and discriminators; deep convolutional GAN \citep{radford2016unsupervised}, progressive GAN \citep{karras2018progressive} and style GAN \citep{karras2019style} are successful ones.
In many real applications, GAN often performs better than the likelihood approach in terms of the quality of generated samples.

In spite of the rapid development of GAN, a theoretical understanding of it remains largely unexplored. Specifically, a generative model typically focuses on providing an estimator only for the generator $\bg$ and does not yield an explicit estimator for the unknown distribution $P_0$. Since the generator is not identifiable, it is crucial to study the convergence rate of the distribution estimator, implicitly defined through $\hat \bg$. This paper aims to bridge this gap by studying the statistical properties of GAN from the viewpoint of estimating a nonparametric distribution. We investigate the convergence rate of a GAN-based estimator for the underlying distribution concentrated around a low-dimensional structure. Through this analysis, our objective is to provide theoretical insights into why GAN outperforms classical nonparametric methods and likelihood approaches in many applications.

Let $Q_\bg$ and $P_{\bg, \sigma}$ denote distributions of $\bg(\bZ)$ and $\bg(\bZ) + \bepsilon$, respectively, where $\bZ \sim P_Z$ and $\bepsilon \sim \cN(\bzero_D, \sigma^2\Id_D)$ are independent.
$Q_\bg$ is often called the pushforward measure of $P_Z$ through the generator $\bg$.
For the data-generating distribution $P_0$, we assume that $P_0 = P_{\bg_0, \sigma_0}$ with a true generator $\bg_0$ and $\sigma_0 \geq 0$. 
We further assume that $\bg_0$ possesses a certain low-dimensional structure and $\sigma_0$ is small enough so that $P_0$ is concentrated around the structure.
This assumption on the true distribution has been thoroughly investigated by \cite{chae2023likelihood}, inspired by recent articles on structured distribution estimation \citep{genovese2012manifold, genovese2012minimax, puchkin2022structure, aamari2019nonasymptotic, divol2020minimax}.
Once the true generator $\bg_0$ possesses a low-dimensional structure that DNN can efficiently capture, deep generative models are highly appropriate for statistical inferences.
We consider a composite structure \citep{horowitz2007rate, juditsky2009nonparametric} which has recently been studied in deep supervised learning \citep{bauer2019deep, schmidt2020nonparametric}.
Then, the corresponding distribution $Q_0 = Q_{\bg_0}$ inherits the structure of $\bg_0$. Details are described further in Section \ref{sec:true}. 
Although the structural assumption on the distribution through the generator is quite natural, to the best of our knowledge, it has not been studied in the literature except for the work presented in \cite{chae2023likelihood}. Similarly to \cite{chae2023likelihood}, we adopt this structural assumption in order to develop a statistical theory that explains the benefits of deep generative models and GANs.

Under the above setting, it would be more reasonable to set $Q_0$, rather than $P_0$, as the target distribution to be estimated because $\bepsilon$ is a noise.
Once an estimator $\hat \bg$ is constructed, one can define an estimator for $Q_0$ as $\hat Q = Q_{\hat\bg}$.
To evaluate the performance of the estimation, we primarily consider the $L^1$-Wasserstein metric.
The metric was originally inspired by the problem of optimal mass transportation \citep{villani2003topics} and has been widely adopted as an evaluation measure in distribution estimation problems \citep{nguyen2013convergence, chae2019bayesianConsistency, wei2022convergence}.
When $\bg_0$ possesses a composite structure with parameters $(t_i, \beta_i)_{i=0}^q$, see Section \ref{sec:true} for details, we construct a GAN-based estimator that achieves the convergence rate 
\bean
    \max_{i \in \{0, \ldots, q\}} n^{-\frac{\beta_i}{2\beta_i + t_i}} + \sigma_0
\eean
up to a logarithmic factor; see Theorem \ref{thm:rate-composition}.
Note that the rate does not explicitly depend on the dimensions $D$ and $d$.
Under the same assumption, \cite{chae2023likelihood} obtained the rate 
\bean
    \max_{i \in \{0, \ldots, q\}} n^{-\frac{\beta_i}{2(\beta_i + t_i)}} + \sigma_0
\eean
up to a logarithmic factor using a VAE-type estimator.
Note that our convergence rate is strictly faster than the rate obtained in \cite{chae2023likelihood}, which is derived using the sharp probability inequality for likelihood ratios developed by \cite{wong1995probability}. Based on this observation, we conjecture that the convergence rate for the VAE-type estimator in \cite{chae2023likelihood} cannot be improved. If this conjecture holds true, our theory will provide valuable insights into the reasons why GAN outperforms VAE.

For the class of structured distributions described above, we also obtain a lower bound 
\bean
    \max_{i \in \{0, \ldots, q\}} n^{-\frac{\beta_i}{2\beta_i + t_i - 2}}
\eean
for the minimax convergence rate; see Theorem \ref{thm:lower-bound}.
When $\sigma_0$ is small enough, this lower bound is only slightly smaller than the rate achieved by a GAN-based estimator.
That is, the convergence rate of a GAN-based estimator obtained in this paper is at least very close to the minimax optimal rate.
As discussed after Theorem \ref{thm:lower-bound}, we conjecture that our lower bound cannot be improved in general, and thus, there might be room for improving the upper bound.

Besides the convergence rate with respect to the $L^1$-Wasserstein distance, we also investigate the convergence rate for $d_{\cF_0}(\hat Q, Q_0)$ with a general integral probability metric $d_{\cF_0}$, as defined in \eqref{eq:ipm}; see Theorem \ref{thm:rate-ipm}. The $L^1$-Wasserstein distance corresponds to a special case where $\cF_0$ is the class of every function with Lipschitz constant bounded by 1. For another example, $\cF_0$ can be chosen as an $\alpha$-Hölder class. Additionally, neural network distances are also natural choices for $d_{\cF_0}$, where the term neural network distance refers to an integral probability metric $d_{\cF_0}$ with $\cF_0$ consisting of neural network functions \citep{arora2017generalization, zhang2018discrimination, bai2019approximability, liu2017approximation}.

It would be worthwhile to highlight several technical novelties of this paper compared to existing theories on GANs. A comprehensive overview of related work can be found in Section \ref{ssec:related-work}.

Firstly, while most existing theories on GANs analyze them from the perspective of nonparametric density estimation, our paper distinguishes itself by focusing on distribution estimation. This allows us to handle both scenarios where the underlying distribution is singular with respect to the Lebesgue measure or possesses a smooth Lebesgue density. In particular, within the framework of existing theory, classical methods such as kernel density estimators and wavelets can achieve the minimax optimal convergence rate. Therefore, their results are insufficient to explain the advantage of GAN compared to classical methods.  In this regard, our theory for GAN is particularly beneficial as it provides a framework that can explain the advantages of using GAN for both density estimation and structured distribution estimation problems. There have been recent articles that explore modifications of classical methods for estimating distributions on manifolds \citep{berenfeld2021density, divol2022measure}. However, it remains unclear whether these methods are suitable for the structured distribution estimation problem addressed in the present paper. The structured distribution estimation considered in our paper involves a substantially richer structure than the manifold structure, as discussed in \cite{chae2023likelihood}.

Another notable technique in the proof of Theorem \ref{thm:rate-composition} lies in the construction of the discriminator class.
In the literature, the function class for the discriminator is identical to the function class defining the evaluation metric.
In case of the $L^1$-Wasserstein, for example, it is the class $\cF_{\rm Lip}$ of every function with Lipschitz constant bounded by one.
In particular, the discriminator class depends solely on the evaluation metric.
On the other hand, the discriminator class in our proof depends not only on the evaluation metric but on the generator architecture.
Although state-of-the-art GAN architectures such as progressive GAN \citep{karras2018progressive} and StyleGAN \citep{karras2019style} are too complicated to render them theoretically tractable, it is crucial for the success of these procedures that discriminator architectures have similar structures to the generator architectures.

In the proof of Theorem \ref{thm:rate-composition}, we carefully construct the discriminator class using the generator class.
In particular, the discriminator class is constructed so that its complexity, expressed through the metric entropy, is of the same order as that of the generator class.
Consequently, the discriminator class becomes a much smaller class than $\cF_{\rm Lip}$, which is the one considered in the literature for obtaining a Wasserstein rate.
By reducing the complexity of the discriminator class, we can significantly improve the convergence rate.

Finally, we would like to mention that once the statement of Theorem \ref{thm:lower-bound} is slightly modified, it might be possible to derive similar lower bounds more easily based on Caffarelli's regularity theory of optimal transport \citep{caffarelli1990interior, urbas1988regularity} and minimax theory for density models \citep{liang2021well, uppal2019nonparametric, niles2022minimax}. More specifically, if $P_Z$ is a uniform distribution on a Euclidean ball in $\bbR^d$ instead of the uniform distribution on the cube $[0,1]^d$ as in Theorem \ref{thm:lower-bound}, Caffarelli's theory provides a useful connection between the density model and generative model, which facilitates an easier proof for the lower bound, see the discussion after Theorem \ref{thm:lower-bound} for more details. However, extending this approach to the case where $P_Z = {\rm Unif}([0,1]^d)$ is not straightforward because the uniform convexity of the support of probability measures involved is a key assumption in Caffarelli's theory. In particular, we may need to construct a sufficiently regular transport map, whose Jacobian determinant is bounded from above and below, from the uniform distribution on a Euclidean ball to the uniform distribution on a cube. Instead of applying the technically involved Caffarelli's regularity theory, we have chosen to directly construct multiple testing based on generators and apply Fano's method to obtain the lower bound. We believe this approach is novel and provides a different perspective.

The remainder of the paper is organized as follows. 
First, we review the literature on the theory of GAN and introduce some notations.
Section \ref{sec:gan} provides a mathematical set-up, including a brief introduction to DNN and GAN.
In Section \ref{sec:true}, we discuss the assumption on the true distribution in depth.
An upper bound for a convergence rate of a GAN-based estimator and a lower bound of minimax convergence rates are investigated in Sections \ref{sec:rate} and \ref{sec:lower-bound}, respectively.
Concluding remarks follow in Section \ref{sec:conclusion}.
All proofs are provided in Supplement.

\subsection{Related statistical theory for GAN} \label{ssec:related-work}

Convergence rates of nonparametric generative models were initially studied in \cite{kundu2014latent} and \cite{pati2011posterior}.
Rather than utilizing DNN, they considered a nonparametric Bayesian approach with a Gaussian process prior on the generator function.

Since the development of GAN by \cite{goodfellow2014generative}, several researchers have studied rates of convergence in deep generative models, particularly focusing on GAN.
An earlier version of \cite{liang2021well} was the first one to study the convergence rate under a GAN framework. More specifically, they considered the Sobolev IPMs to evaluate the estimation performance.
A similar theory has been developed by \cite{singh2018nonparametric}, which was later generalized by \cite{uppal2019nonparametric} using Besov IPMs.
Slightly weaker results were obtained by \cite{chen2020statistical}. 
Although their convergence rates are strictly slower than the minimax optimal rate, they explicitly considered DNN architectures for the generator and discriminator classes.
Convergence rates of the vanilla GAN with respect to the Jensen--Shannon divergence have recently been obtained by \cite{belomestny2021rates}.

These works utilized the framework of nonparametric density estimation to understand GAN.
They evaluated the performance of GAN using integral probability metrics, while classical approaches such as the kernel density estimation focused on other metrics such as the total variation, Hellinger and uniform metrics.
Since the total variation can be viewed as an IPM, some results in the above papers are comparable with that of the classical methods.
In these comparable cases, both approaches achieve the same minimax optimal rate; hence these theories on GAN cannot explain why deep generative models outperform classical nonparametric methods.

\cite{schreuder2021statistical} considered generative models where the target distribution may not possess a Lebesgue density. They assumed that the true distribution is the convolution of $Q_{\bg_0}$ and a general noise distribution for some function $\bg_0: [0,1]^d \to [0,1]^D$. While this assumption is similar to ours, it does not explicitly incorporate the smoothness and composite structure of $\bg_0$. As a result, their result only guarantees that GAN achieves the same rate as the empirical measure. More recently, \cite{tang2023minimax} obtained the minimax rate of distribution estimation under a submanifold assumption using a mixture of GANs. However, as mentioned earlier, the composite structure imposed through the generator function in our paper involves a substantially richer structure than just a manifold structure considered in \cite{tang2023minimax}. For instance, the dimension $t_*$ corresponding to the worst-case component of a composite function (as defined in \eqref{eq:intrinsic-def}) can be much smaller than the dimension of the manifold on which $Q_0$ is supported.

\subsection{Notations}

Maximum and minimum of two real numbers $a$ and $b$ are denoted as $a \vee b$ and $a \wedge b$, respectively.
For $1 \leq p < \infty$, $|\cdot |_p$ denotes the $\ell^p$-norm.
For a real-valued function $f$ and a probability measure $P$, let $P f = \int f(\bx) dP(\bx)$.
$\E$ denotes the expectation when the underlying probability is obvious.
Convolution of two probability measures $P$ and $Q$ are denoted $P*Q$.
We write $c = c(A_1, \ldots, A_k)$ when $c$ depends only on $A_1, \ldots, A_k$.
Uppercase letters, such as $P$ and $\hat P$, refer to probability measures corresponding to densities denoted by their lowercase counterparts: \ie\ $p$ and $\hat p$, respectively.
We write $a \lesssim b$ if $a$ is less than $b$ up to a constant multiplication, where the constant is universal or at least contextually unimportant.
Lastly, $a \asymp b$ indicates $a \lesssim b$ and $b \lesssim a$.

\section{Generative adversarial networks} \label{sec:gan}

For a given class $\cF$ of functions from $\bbR^D$ to $\bbR$, the $\cF$-IPM (\cite{muller1997integral}) between two probability measures $P_1$ and $P_2$ is defined as
\be \label{eq:ipm}
  d_\cF(P_1, P_2) = \sup_{f \in \cF} |P_1 f - P_2 f|.
\ee
For example, if $\cF = \cF_{\rm Lip}$, the class of every function $f: \bbR^D \to \bbR$ satisfying $| f(\bx) - f(\by)| \leq |\bx - \by|_2$ for all $\bx, \by \in \bbR^D$, then the corresponding IPM is the $L^1$-Wasserstein distance by the Kantorovich--Rubinstein duality theorem; see Theorem 1.14 from \cite{villani2003topics}.
H\"{o}lder, or more generally Besov, IPMs have been considered in recent articles for evaluating the performance of the density or distribution estimation; see 
\cite{liang2021well},
\cite{uppal2019nonparametric},
\cite{singh2018nonparametric} and
\cite{tang2023minimax}.

Let $\cG$ be a class of functions from $\cZ\subset \bbR^d$ to $\bbR^D$, and $\cF$ be a class of functions from $\bbR^D$ to $\bbR$.
Two classes $\cG$ and $\cF$ are referred to as the \emph{generator} and \emph{discriminator} classes, respectively.
For given discriminator and generator classes, we define a GAN-based estimator $\hat\bg$ as the minimizer of $d_\cF (Q_\bg, \bbP_n)$ over $\cG$, where $\bbP_n$ is the empirical measure based on the $D$-dimensional observations $\bX_1, \ldots, \bX_n$.
That is, the estimator $\hat\bg \in \cG$ is such that
\be\label{eq:GAN-estimator}
	d_\cF(Q_{\hat\bg}, \bbP_n) \leq \inf_{\bg \in \cG} d_\cF (Q_\bg, \bbP_n) + \epsilon_{\rm opt}.
\ee
Here, the optimization error $\epsilon_{\rm opt} \geq 0$ is a prespecified number.
An estimator satisfying \eqref{eq:GAN-estimator} is of our primary interest.
Although the vanilla GAN \citep{goodfellow2014generative} is not of the form \eqref{eq:GAN-estimator}, the formulation \eqref{eq:GAN-estimator} is quite general to include various GANs popularly used in practice \citep{arjovsky2017wasserstein, li2017mmd, mroueh2017sobolev}.
At the population level, \eqref{eq:GAN-estimator} can be viewed as a method to solve the following minimization
\bean
  \minimize_{\bg \in \cG} d_\cF(Q_\bg, P_0)
\eean
because one may expect $\E d_\cF(\bbP_n, P_0) \to 0$, where the expectation is taken with respect to $P_0$.
Since the convergence rate for $\E d_\cF(\bbP_n, P_0)$ might be very slow, however, a careful analysis is necessary.
In particular, in Section \ref{sec:rate}, we will separate the evaluation metric from the $\cF$-IPM, the one defined through the discriminator class.

In practice, both the generator and discriminator classes are modeled using deep neural networks.
To be specific, let $\rho(x) = x \vee 0$ be the ReLU activation function \citep{glorot2011deep}.
We focus on the ReLU in this paper, but other activation functions can also be used as long as a suitable approximation property holds \citep{ohn2019smooth}.
For vectors $\bv = (v_1, \ldots, v_r)$ and $\bz = (z_1, \ldots, z_r)$, define $\rho_\bv(\bz) = (\rho(z_1-v_1), \ldots, \rho(z_r - v_r))$.
For a nonnegative integer $L$ and $\bp = (p_0, \ldots, p_{L+1}) \in \bbN^{L+2}$, a neural network function with the network architecture $(L, \bp)$ is any function $\bff: \bbR^{p_0} \to \bbR^{p_{L+1}}$ such that
\be\label{eq:dnn}
	\bz \mapsto \bff(\bz) = W_L \rho_{\bv_L} W_{L-1} \rho_{\bv_{L-1}} \cdots W_1 \rho_{\bv_1} W_0 \bz,
\ee
where $W_i \in \bbR^{p_{i+1} \times p_i}$ and $\bv_i \in \bbR^{p_i}$.
Let $\cD(L, \bp, s, F)$ be the collection $\bff$ from \eqref{eq:dnn} satisfying
\bean
	\max_{j=0, \ldots,L} |W_j|_\infty \vee |\bv_j|_\infty \leq 1,
	\quad	
	\sum_{j=1}^L |W_j|_0 + |\bv_j|_0 \leq s
	\quad \text{and} ~~ \| \bff \|_\infty \leq F,
\eean
where $|W_j|_\infty$ and $|W_j|_0$ denote the maximum-entry norm and the number of nonzero elements of the matrix $W_j$, respectively, and $\|\bff\|_\infty = \| |\bff(\bz)|_\infty \|_\infty = \sup_\bz |\bff(\bz)|_\infty$.

When the generator class $\cG$ consists of neural network functions, we call the corresponding class $\cQ = \{Q_\bg: \bg \in \cG\}$ as a \emph{deep generative model}.
In this sense, GAN can be viewed as a method for estimating the parameters in deep generative models.
In the literature concerning the variational autoencoder, a collection of $P_{\bg, \sigma}$ is often called a deep generative model as well.

\section{Assumptions on the true distribution} \label{sec:true}

In this section, we address assumptions on the true distribution $P_0$.
As mentioned in the introduction, we assume that $P_0 = P_{\bg_0, \sigma_0}$ for some function $\bg_0: \cZ \to \bbR^D$ and $\sigma_0 \geq 0$.
Furthermore, we assume that $\bg_0$ possesses a structure that DNN can efficiently capture.
As long as $\sigma_0$ is not too large, the true distribution $P_0$ inherits the structure of $\bg_0$, which enables efficient estimation of it (or $Q_0 = Q_{\bg_0}$).
Note that it is much more convenient to impose a structure on the generator rather than directly on the density function because there is no constraint on the functional form of the generator.

We suppose that $\bg_0$ belongs to a class of structured functions.
More specifically, we consider a class of composite functions for which deep generative models have benefits.
For positive numbers $\beta$ and $K$, let $\cH^\beta_K(A)$ be a class of all functions from $A$ to $\bbR$ with $\beta$-\Holder norm bounded by $K$. See \cite{van1996weak} and \cite{gine2016mathematical} for the definition of \Holder space.
Consider a function $\bg: \bbR^d \to \bbR^D$ as follows:
\be \label{eq:composition}
	\bg = \bh_q \circ \bh_{q-1} \circ \cdots \circ \bh_1 \circ \bh_0
\ee
with $\bh_i: (a_i, b_i)^{d_i} \to (a_{i+1}, b_{i+1})^{d_{i+1}}$ and $\bh_i = (h_{i1}, \ldots, h_{id_{i+1}})$. Here, $d_0=d$ and $d_{q+1}=D$.
Let $t_i$ be the maximal number of variables on which each of the $h_{ij}$ depends.
Let $\cG_0(q, \bd, \bt, \bbeta, K)$ be a collection of functions of the form \eqref{eq:composition} satisfying $h_{ij} \in \cH^{\beta_i}_K\big( (a_i, b_i)^{t_i} \big)$ and $|a_i| \vee |b_i| \leq K$, where $\bd = (d_0, \ldots, d_{q+1})$, $\bt = (t_0, \ldots, t_q)$ and $\bbeta = (\beta_0, \ldots, \beta_q)$.
Let
\be \label{eq:intrinsic-def}
	\widetilde \beta_i = \beta_i \prod_{l=i+1}^q (\beta_l \wedge 1), 
	\quad
	i_* = \argmax_{i \in \{0, \ldots,q\}} \frac{t_i}{\widetilde \beta_i},
	\quad
	\beta_* = \widetilde\beta_{i_*}
	\quad
	\text{and} ~~ t_* = t_{i_*}.
\ee
Note that the decomposition of the form \eqref{eq:composition} for a given function $\bg$ may not be unique. For example, the composite function $\bh_q \circ \bh_{q-1} \circ \cdots \circ \bh_1 \circ \bh_0$ can be seen as a single function $\widetilde \bh_0 = (\widetilde h_{01}, \ldots, \widetilde h_{0D})$ with $\widetilde h_{0j} \in \cH_{K'}^{\min\{\beta_0, \ldots, \beta_q\}} ((a_0, b_0)^{t_0})$ for a large enough constant $K'$. Also, there might be several maximizers for the map $i \mapsto t_i / \widetilde\beta_i$. In this case, $i_*$ can be defined as any maximizer.

The class $\cG_0 = \cG_0(q, \bd, \bt, \bbeta, K)$ has been extensively studied in recent articles on deep supervised learning to demonstrate the benefits of DNN in estimating a nonparametric function \citep{schmidt2020nonparametric, bauer2019deep}.
As studied in \cite{chae2023likelihood}, a composite structure can naturally be translated to unsupervised learning problems through the distribution class $\cQ_0 = \{Q_\bg: \bg \in \cG_0\}$.
For example, when $d=D$ and $\cG_0$ consists of functions of the form $\bg(\bz) = (g_1(z_1), \cdots, g_d(z_d))$, where $\bz = (z_1, \ldots, z_d)$ and $g_j$, $j = 1, \ldots, d$, is a univariate function, then $t_* = 1$ and the corresponding $\cQ_0$ becomes a class of product distributions.
If $\cG_0$ consists of $\beta$-\Holder functions with $\beta > 1$, the support of $P_Z$ is uniformly convex and its density is bounded from above and below, the corresponding $\cQ_0$ contains distributions possessing a strictly positive $(\beta-1)$-\Holder density on a bounded and uniformly convex subset of $\bbR^D$.
This fact is based on the well-established regularity theory of optimal transport; see Theorem 12.50 of \cite{villani2008optimal} for details. It is important to note that the uniform convexity assumption cannot be relaxed within Caffarelli's regularity theory.

In the literature, structured distribution estimation has predominantly been studied within the framework of manifold structure \citep{tang2023minimax, puchkin2022structure}. However, the composite structure introduced through the generator function in our approach incorporates various interesting low-dimensional structures that are not captured by the manifold structure alone. For instance, consider the example mentioned earlier, where $\bg(\bz) = (g_1(z_1), \cdots, g_d(z_d))$. In this case, the dimension $t_*$, defined as the worst-case component of the composite function (as in \eqref{eq:intrinsic-def}), can be much smaller than the dimension $d$ of the manifold on which $\bg(\bZ)$ is supported. This highlights the richer and more flexible structure captured by the composite approach compared to the  manifold structure.

If $d < D$ and $\bg_0$ is sufficiently smooth, the distribution $Q_0$ is singular with respect to the Lebesgue measure on $\mathbb{R}^D$. However, the distribution $P_0$ possesses a Lebesgue density provided that $\sigma_0 > 0$. We would like to emphasize that the main theorems in Section \ref{sec:rate} hold for all values of $\sigma_0$ in the interval $[0,1]$. With regard to the noise level $\sigma_0$, it would be worthwhile to discuss two different regimes, as also discussed in Section 3.6 of \cite{chae2023likelihood}.

Firstly, consider the case where $\sigma_0$ is a fixed positive constant. In this case, our results do not provide a meaningful convergence rate. The problem of estimating $Q_0$ with additive noise is commonly referred to as the deconvolution problem, and it has been extensively studied in the literature \citep{fan1991optimal, alexander2009deconvolution, genovese2012manifold, nguyen2013convergence} under the assumption of fixed $\sigma_0$. It is worth noting that the estimation problem in this setting is intrinsically very difficult, and this difficulty is often expressed mathematically through the logarithmic minimax rates. While a GAN-based estimator might achieve such a logarithmic convergence rate, we do not pursue its study in the present paper, as our primary focus is on the regime where $\sigma_0$ is small. In particular, we believe that a theory with such a slow convergence rate would not be suitable for explaining the amazing performance of deep generative models.

A small $\sigma_0$ regime can be mathematically expressed as $\sigma_0 \to 0$ with a suitable rate as $n \to \infty$. In this regime, it is possible to obtain a fast convergence rate for estimating $Q_0$, as guaranteed by our theory. While the data-generating distribution $P_0$ depends on the sample size $n$, the theorems in Section \ref{sec:rate} hold for all $n$, ensuring clear interpretation of the results. It is worth noting that such sample size-dependent true distributions have been extensively studied in modern high-dimensional statistics \citep{buhlmann2011statistics, wainwright2019high}, and our setup can be understood within similar contexts. Although our setup and estimation problems may differ slightly, there have been several recent articles that assume data are concentrated around a small neighborhood of a manifold, and these neighborhoods shrink to the manifold as the sample size increases; see 
\cite{puchkin2022structure},
\cite{aamari2019nonasymptotic},
\cite{aamari2018stability},
\cite{divol2021minimax},
\cite{jiao2023deep} and
\cite{berenfeld2022estimating}
for relevant discussions in this direction.

\section{Convergence rate of a GAN-based estimator} \label{sec:rate}

Although a strict minimization of the map $\bg \mapsto d_\cF(Q_\bg, \bbP_n)$ is computationally intractable, several heuristic approaches are available to approximate the solution to \eqref{eq:GAN-estimator}.
In this section, we investigate the convergence rate of $\hat Q =  Q_{\hat\bg}$ under the assumption that the computation of it is possible.
A goal is to find a sharp upper bound for $\E d_{\rm eval} (\hat Q, Q_0)$, where $d_{\rm eval}$ is a metric evaluating the performance of the estimation.
It is desirable that the rate is independent of $D$ and $d$, and it adapts to the structure of $P_0$.
We consider an arbitrary evaluation metric $d_{\rm eval}$ for generality; the $L^1$-Wasserstein distance is of primary interest.
Recall that the $L^p$-Wasserstein distance between two Borel probability measures $P$ and $Q$ on $\bbR^D$ is defined as
\bean
  W_p(P, Q) = \inf_\pi \left( \int |\bx - \by|_2^p d\pi(\bx, \by)\right)^{1/p}
\eean
for $p \geq 1$, where the infimum is taken over every coupling $\pi$ of $P$ and $Q$.
As mentioned earlier, $W_1$ allows the IPM representation $W_1(P, Q) = d_{\cF_{\rm Lip}}(P, Q)$ by the Kantorovich--Rubinstein duality, which makes it convenient to utilize $W_1$ as an evaluation metric in a GAN framework.
Although a more general duality theorem is well-known for $p > 1$ (\cite{villani2003topics}, Theorem 1.3), the IPM representation of $W_p$ is available only for $p=1$.

In literature, the evaluation metric is often identified with $d_\cF$, the IPM defined through the discriminator class $\cF$.
In this sense, when $d_{\rm eval}$ is the $L^1$-Wasserstein metric $W_1$, a natural candidate for the discriminator class $\cF$ might be $\cF_{\rm Lip}$, the class of those functions  whose Lipschitz constant is bounded by 1.
Indeed, it is the original motivation of the Wasserstein GAN to minimize the objective function
\be\label{eq:wgan}
  W_1(Q_\bg, \bbP_n) = d_{\cF_{\rm Lip}}(Q_\bg, \bbP_n)
\ee
over $\bg \in \cG$.
In the original article of the Wasserstein GAN \citep{arjovsky2017wasserstein}, \eqref{eq:wgan} is minimized after replacing $\cF_{\rm Lip}$ by a class $\cF$ of neural network functions.
The replacement was only for computational tractability.
Although minimizing the map $\bg \mapsto d_\cF (Q_\bg, \bbP_n)$ with a neural network class $\cF$ is still challenging, several heuristic approaches can be employed to approximate the solution.

Suppose that computing the minimizer of \eqref{eq:wgan}, say $\hat Q^W$, is possible.
It is natural to ask whether $\hat Q^W$ is a decent estimator, theoretically at least.
If the generator class $\cG$ is large enough, for example, $\hat Q^W$ is expected to be close to the empirical measure.
Consequently, the convergence rates of $\hat Q^W$ and $\bbP_n$ would be the same.
\cite{schreuder2021statistical} utilized this idea to prove that $\hat Q^W$ performs at least as good as $\bbP_n$ does.
The convergence rate of the empirical measure with respect to the Wasserstein metric is well-known in the literature. \cite{fournier2015rate} have shown that
\be\label{eq:rate-empirical}
  \E W_1(\bbP_n, P_0) \lesssim
    \left\{\begin{array}{ll}
      n^{-1/2} & ~~\text{if $D=1$}
      \\
      n^{-1/2} \log n & ~~\text{if $D=2$}
      \\
      n^{-1/D} & ~~\text{if $D > 2$.}
    \end{array}\right.
\ee
See also \cite{weed2019sharp}.
The rate \eqref{eq:rate-empirical} becomes slower as $D$ increases, suffering from the curse of dimensionality.
Although $\bbP_n$ adapts to a certain intrinsic dimension and achieves the minimax rate in some sense \citep{weed2019sharp, singh2018minimax}, it is possible to construct a more efficient estimator, particularly when the underlying distribution possesses a smooth structure.

$\hat Q^W$ may perform better than the empirical measure if $\cG$ is not too large, but the theoretical analysis of this would be quite challenging.
Furthermore, practical estimators might be fundamentally different from $\hat Q^W$ because it does not take crucial features of state-of-the-art methods into account.
For successful GAN approaches, for example, the structures of the generator and discriminator architectures are closely related. In particular, the complexities of the two architectures are similar.
On the other hand, the discriminator class  $\cF=\cF_{\rm Lip}$, used in the construction of $\hat Q^W$, has no connection with the generator class.
In this sense, it is difficult to view $\hat Q^W$ as a suitable estimator to be theoretically analyzed.

In conclusion, $d_{\rm eval}$ is not necessarily identical to $d_\cF$ in our analysis.
Nonetheless, $d_\cF$ should be close to $d_{\rm eval}$ in some sense because GAN constructs an estimator by minimizing $d_\cF(Q_\bg, \bbP_n)$ over $\bg \in \cG$.
This is specified as condition (iv) of Theorem \ref{thm:rate-general}: $d_\cF$ needs to be close to $d_{\rm eval}$ only on a relatively small class of distributions.

\begin{theorem}\label{thm:rate-general}
Suppose that $\bX_1, \ldots, \bX_n$ are \iid\ random vectors following $P_0 = Q_0 * \cN(\bzero_D, \sigma_0^2 \Id_D)$ for some distribution $Q_0$ (not necessarily of the form $Q_{\bg_0}$) and $\sigma_0 \geq 0$.
For given generator class $\cG$, discriminator class $\cF$ and an estimator $\hat Q = Q_{\hat\bg}$ with $\hat\bg \in \cG$, suppose that
\be\begin{split} \label{eq:general-condition}
	{\rm (i)}~ & \inf_{\bg\in\cG} d_{\rm eval} (Q_\bg, Q_0) \leq \epsilon_1
	\\
	{\rm (ii)}~ & d_\cF(\hat Q, \bbP_n) \leq \inf_{\bg \in \cG} d_\cF (Q_\bg, \bbP_n) + \epsilon_2
	\\
	{\rm (iii)}~ & \E d_\cF(\bbP_n, P_0) \leq \epsilon_3
	\\
	{\rm (iv)}~ & |d_{\rm eval}(Q_1, Q_2) - d_\cF(Q_1, Q_2)| \leq \epsilon _4 
	\quad \text{$\forall Q_1, Q_2 \in \cQ \cup \{Q_0\}$,}
\end{split}\ee
where $\cQ = \{Q_\bg: \bg \in \cG\}$ and $\epsilon_j \geq 0$.
Then,
\bean
	\E d_{\rm eval}(\hat Q, Q_0) \leq 2 d_\cF(P_0, Q_0) + 5\epsilon_1 + \epsilon_2 + 2\epsilon_3 + 3\epsilon_4.
\eean
\end{theorem}

Note that similar inequalities to the statement of Theorem \ref{thm:rate-general} have been explicitly or implicitly considered in the literature to analyze the theoretical properties of GANs \citep{chen2020statistical, biau2021some, schreuder2021statistical, belomestny2021rates}. Theorem \ref{thm:rate-general} is a slight modification of these existing results, with the modification favoring our analysis. The proof of Theorem \ref{thm:rate-general} does not significantly differ from the proofs in the literature.

Two quantities $\epsilon_1$ and $\epsilon_3$ are closely related to the complexity of $\cG$ and $\cF$, respectively.
In particular, $\epsilon_1$ represents an error for approximating $Q_0$ by distributions of the form $Q_\bg$ over $\bg \in \cG$. The larger the generator class $\cG$ is, the smaller the approximation error is
\citep{yarotsky2017error, telgarsky2016benefits, ohn2019smooth}.
Similarly, $\epsilon_3$ gets larger as the complexity of $\cF$ increases.
Techniques for bounding $\E d_\cF(\bbP_n, P_0)$ are well-known in empirical process theory \citep{van1996weak, gine2016mathematical}.
The second error term $\epsilon_2$ is nothing but the optimization error.
The fourth term $\epsilon_4$ is the deviance between the evaluation metric $d_{\rm eval}$ and $\cF$-IPM over $\cQ \cup \{Q_0\}$, connecting $d_\cF$ and $d_{\rm eval}$.

Finally, the term $d_\cF(P_0, Q_0)$ in the assertion of Theorem \ref{thm:rate-general} depends primarily on $\sigma_0$.
If $\cF \subset \cF_{\rm Lip}$, for example, one can easily prove that 
\be\label{eq:sigma-ineq}
	d_\cF(P_0, Q_0) \leq W_1(P_0, Q_0) \leq W_2(P_0, Q_0) \lesssim \sigma_0,
\ee
where the third inequality holds by well-known formula \citep{givens1984class}. As another example, if $\cF$ consists of twice continuously differentiable functions with suitably bounded derivatives, we have
\be\begin{split}\label{eq:sigma2}
	& |P_0 f - Q_0 f| 
	= \big| \E [ f(\bY + \bepsilon) - f(\bY) ] \big|
	\\
	&\approx \Big| \E \Big[ \bepsilon^T \nabla f(\bY) + \frac{1}{2} \bepsilon^T \nabla^2 f(\bY) \bepsilon \Big] \Big|
	\asymp \sigma_0^2,
\end{split}\ee
for $f \in \cF$ where $\bY \sim Q_0$ and $\bepsilon \sim \cN(\bzero_D, \sigma_0^2 \Id)$ are independent random vectors.
Hence, $d_\cF (P_0, Q_0) \lesssim \sigma_0^2$, which gives a better bound than \eqref{eq:sigma-ineq} for a small enough $\sigma_0$.

Ignoring the optimization error, suppose for a moment that $\cG$ is given and we need to choose a suitable discriminator class to minimize $\epsilon_3 + \epsilon_4$ in Theorem \ref{thm:rate-general}.
We focus on the case of $d_{\rm eval} = W_1$.
One can easily make $\epsilon_4 = 0$ by taking $\cF = \cF_{\rm Lip}$.
In this case, however, $\epsilon_3$ would be too large because $\E W_1 (\bbP_n, P_0) \asymp n^{-1/D}$ for $D > 2$; \cf\ Eq.\ \eqref{eq:rate-empirical}.
That is, $\cF_{\rm Lip}$ might be too large to be used as a discriminator class. The discriminator class $\cF$ should be much smaller than $\cF_{\rm Lip}$ to obtain a fast convergence rate.
To achieve this goal, we construct $\cF$ so that both $\epsilon_3$ and $\epsilon_4$ are small enough.
Such discriminator class can be constructed as, for example,
\be\label{eq:discriminator}
	\cF = \Big\{f_{Q_1, Q_2}: Q_1, Q_2 \in \cQ \cup \{Q_0\} \Big\},
\ee
where $f_{Q_1, Q_2}$ is a (approximate) maximizer of $|Q_1 f - Q_2 f|$ over $f \in \cF_{\rm Lip}$.
In this case, $\epsilon_4$ vanishes and the convergence rate of $\hat Q$ will be determined solely by $\epsilon_1$, $\epsilon_3$ and $\sigma_0$.
Furthermore, the complexity of $\cF$ would roughly be the same as that of $\cG\times\cG$.
If the complexity of a function class is expressed through a metric entropy, the logarithmic covering number, the complexities of $\cG$ and $\cF$ are of the same order.
In this case, three quantities $\epsilon_1$, $\epsilon_3$ and $\sigma_0$ can roughly be interpreted as the approximation error, estimation error and noise level, respectively.
While we cannot control the noise level $\sigma_0$, both the approximation and estimation errors depend on the complexity of $\cG$. Hence, a suitable choice of it is important to achieve a fast convergence rate.

As described in Section \ref{sec:true}, we consider a class of composite functions for the true generator.
Let 
\bean
	\cQ_0 = \cQ_0(q, \bd, \bt, \bbeta, K) = \Big\{Q_\bg: \bg \in \cG_0(q, \bd, \bt, \bbeta, K) \Big\},
\eean
where $\cG_0 = \cG_0(q, \bd, \bt, \bbeta, K)$ is defined as in Section \ref{sec:true}.
Although not strictly necessary, it would be convenient to regard quantities $(q, \bd, \bt, \bbeta, K)$ as constants independent of $n$.
In the forthcoming Theorem \ref{thm:rate-composition}, we obtain a Wasserstein convergence rate of $\hat Q$ under the assumption that $Q_0 \in \cQ_0$.

\begin{theorem}\label{thm:rate-composition}
Suppose that $\bX_1, \ldots, \bX_n$ are \iid\ random vectors following $P_0 = P_{\bg_0, \sigma_0}$ for some $\sigma_0 \leq 1$ and $\bg_0 \in \cG_0(q, \bd, \bt, \bbeta, K)$.
Then, there exist a generator class $\cG = \cD(L, \bp, s, K\vee 1)$ and a discriminator class $\cF \subset \cF_{\rm Lip}$ such that for an estimator $\hat Q$ satisfying \eqref{eq:GAN-estimator}, 
\normalsize
\be\begin{split}\label{eq:rate-composition}
	\sup_{Q_0 \in \cQ_0} \E W_1(\hat Q, Q_0)
	\leq C \bigg\{ n^{-\frac{\beta_*}{2\beta_* + t_*}}  ( \log n )^{\frac{3\beta_*}{2 \beta_* + t_*}} 
	+ \sigma_0 + \epsilon_{\rm opt} \bigg\},
\end{split}\ee
where $C = C(q, \bd, \bt, \bbeta, K)$.
\end{theorem}

Theorem \ref{thm:rate-composition} only considers a Gaussian additive noise, but the assumption $P_0 = P_{\bg_0, \sigma_0} = Q_{\bg_0} * \cN(\bzero_D, \sigma_0^2 \Id_D)$ can be relaxed in various ways. In the proof of Theorem \ref{thm:rate-composition}, with regard to the data distribution $P_0$, we only need a bound $d_\cF(P_0, Q_0) \lesssim \sigma_0$ as in \eqref{eq:sigma-ineq} and that $f(\bX_i)$ is a sub-Gaussian variable for every $f \in \cF_{\rm Lip}$, with $f(\bzero_D) = 0$, where the sub-Gaussian parameter $\sigma$ is independent of $f$. Therefore, for example, the normal distribution $\cN(\bzero_D, \sigma_0^2 \Id_D)$ for the noise distribution can be replaced by any sub-Gaussian distribution with variance $\sigma_0^2$.

In Theorem \ref{thm:rate-composition}, both the generator class $\cG$ and discriminator class $\cF$ depend solely on the sample size $n$ and the parameters $(q, \bd, \bt, \bbeta, K)$, independent of $Q_0$ or $P_0$. Moreover, from the proof, it can be deduced that the network parameters $(L, \bp, s)$ of the generator class can be chosen such that $L \lesssim \log n$, $|\bp|_\infty \lesssim n^{t_*/(2\beta_* + t_*)}$ and $s \lesssim n^{t_*/(2\beta_* + t_*)} \log n$, where the constants in $\lesssim$ depend only on $(q, \bd, \bt, \bbeta, K)$.

Ignoring the optimization error $\epsilon_{\rm opt}$, the rate \eqref{eq:rate-composition} consists of the two terms, $\sigma_0$ and $n^{-\beta_*/(2\beta_* + t_*)}$ up to a logarithmic factor.
If $\sigma_0 \lesssim n^{-\beta_*/(2\beta_* + t_*)}$, it can be absorbed into the polynomial term. Therefore, $\hat Q$ achieves the rate $n^{-\beta_*/(2\beta_* + t_*)}$ when $\sigma_0$ is small enough.
Note that this rate often appears in nonparametric smooth function estimation, balancing the approximation and estimation errors.

Under the condition given in Theorem \ref{thm:rate-composition}, \cite{chae2023likelihood} considered a likelihood approach to study the benefit of the deep generative model.
More specifically, they obtained a convergence rate of a sieve MLE based on a Gaussian mixture density $p_{\bg, \sigma}$.
Note that the density $p_{\bg, \sigma}$ is concentrated around a small neighborhood of a low-dimensional structure induced by $\bg$. As a result, likelihood approaches might be highly unstable due to the singularity issue.
To overcome this problem, \cite{chae2023likelihood} considered a sieve MLE based on the perturbed data $\widetilde\bX_i = \bX_i + \widetilde\bepsilon_i$, where $\widetilde\bepsilon_i$ is an artificial noise alleviating the problem caused by the singularity.
They proved that a sieve MLE with an optimal perturbation achieves the Wasserstein rate $n^{-\beta_*/2(\beta_* + t_*)} + \sigma_0$.
The rate was obtained based on the perturbed data, thus it was conjectured to be suboptimal.
Theorem \ref{thm:rate-composition} shows that GAN achieves a strictly faster convergence rate than $n^{-\beta_*/2(\beta_* + t_*)} + \sigma_0$, the one obtained by \cite{chae2023likelihood}.
Hence, the rate of a likelihood approach in \cite{chae2023likelihood} is sub-optimal.

In some special cases, the convergence rate of a GAN-based estimator obtained from Theorem \ref{thm:rate-composition} can be strictly worse than the rate achieved by the empirical measure. For instance, when $\beta_* = 1$, $D = d = t_* > 2$, and $\sigma_0 = \epsilon_{\rm opt} = 0$, the rate \eqref{eq:rate-composition} simplifies to $n^{-1/(d+2)}$ (up to a logarithmic factor), while the empirical measure $\bbP_n$ achieves a strictly faster rate of $n^{-1/d}$, as shown in \eqref{eq:rate-empirical}. However, by choosing different generator and discriminator classes, it is possible to obtain a GAN-based estimator with a convergence rate equal to that of the empirical measure. Specifically, we can apply Theorem \ref{thm:rate-general} with $\cF = \cF_{\rm Lip}$. In this case, $\epsilon_4 = 0$ since $d_{\rm eval} = W_1$. Moreover, by selecting a large enough $\cG$, \ie\, increasing the depth, width, and number of nonzero parameters, we can make $\epsilon_1$ arbitrarily small. As a result, Theorem \ref{thm:rate-general} and Eq.\ \eqref{eq:sigma-ineq} yield the bound $\E W_1(\hat Q, Q_0) \lesssim \sigma_0 + \epsilon_{\rm opt} + \E W_1(\bbP_n, P_0)$. In summary, if 
\be \label{eq:two-case}
	\E W_1 (\bbP_n, P_0) \leq n^{-\frac{\beta_*}{2\beta_* + t_*}}  ( \log n )^{\frac{3\beta_*}{2 \beta_* + t_*}},
\ee
choosing alternative generator and discriminator classes results in an estimator with a convergence rate better than that in Theorem \ref{thm:rate-composition}. We consider this alternative choice in the statement of Theorem \ref{thm:rate-ipm}.

So far, we have focused on the case $d_{\rm eval} = W_1$.
In the remainder of this section, we consider a general IPM as an evaluation metric.
The function space defining the evaluation metric will be denoted $\cF_0$, hence $d_{\rm eval} = d_{\cF_0}$.

\begin{theorem} \label{thm:rate-ipm}
Suppose that $\bX_1, \ldots, \bX_n$ are \iid\ random vectors following $P_0 = P_{\bg_0, \sigma_0}$ for some $\bg_0 \in \cG_0(q, \bd, \bt, \bbeta, K)$ and $\sigma_0 \leq 1$.
Let $\cF_0$ be a class of Lipschitz continuous functions from $\bbR^D$ to $\bbR$ with Lipschitz constant bounded by a constant $C_1 > 0$.
Then, there exist a generator class $\cG = \cD(L, \bp, s, K\vee 1)$ and a discriminator class $\cF$ such that $\hat Q$ defined as in \eqref{eq:GAN-estimator} satisfies
\normalsize
\be\begin{split} \label{eq:rate-ipm}
	\sup_{Q_0 \in \cQ_0} \E d_{\cF_0}(\hat Q, Q_0)	
	\leq C_2 \bigg\{ \sigma_0 + \epsilon_{\rm opt}
	+ n^{-\frac{\beta_*}{2\beta_* + t_*}}  ( \log n )^{\frac{3\beta_*}{2 \beta_* + t_*}} \bigg\},
\end{split}\ee
where $C_2 = C_2(q, \bd, \bt, \bbeta, K, C_1)$.
Alternatively, if we take $\cF = \cF_0$ and the depth, width and number of nonzero parameters of $\cG$ are large enough, the estimator satisfies
\be \label{eq:rate-ipm2}
    \E d_{\cF_0}(\hat Q, Q_0) \leq C_3\Big\{ \sigma_0 + \epsilon_{\rm opt} + \E d_{\cF_0} (\bbP_n, P_0) \Big\},
\ee
where $C_3 = C_3(D, C_1)$.
\end{theorem}

Note that the rate in \eqref{eq:rate-ipm} is slower than that in \eqref{eq:rate-ipm2} if
\bean
	\E d_{\cF_0} (\bbP_n, P_0) \lesssim n^{-\frac{\beta_*}{2\beta_* + t_*}}  ( \log n )^{\frac{3\beta_*}{2 \beta_* + t_*}}.
\eean
It is unclear whether it is possible to construct an estimator which achieves the rate of the minimum of two rates in \eqref{eq:rate-ipm} and \eqref{eq:rate-ipm2}.
This interesting problem is left as a topic for future research.

When $\cF_0$ consists of neural networks, $\cF_0$-IPM is often called a \emph{neural network distance}.
Although it is not a standard choice, neural network distances can serve as an evaluation metric.
In particular, convergence in a neural network distance guarantees a weak convergence under mild assumptions \citep{zhang2018discrimination}.
If $\cF_0 = \cD(L_0, \bp_0, s_0, \infty)$, then it is not difficult to see that $\E d_{\cF_0} (\bbP_n, P_0) \lesssim \sqrt{s_0/n}$ up to a logarithmic factor.
This can be proved using a well-known empirical process theory combined with metric entropy of deep neural networks (See Lemma 5 from \cite{schmidt2020nonparametric}).
Therefore, if $s_0 \gg n^{t_* / (2\beta_* + t_*)}$, the right hand side of \eqref{eq:rate-ipm} provides a strictly faster rate than \eqref{eq:rate-ipm2}.

Another important class of metrics is a \Holder IPM.
When $\cF_0 = \cH^\alpha_1([-K, K]^D)$ for some $\alpha > 0$, it is well-known that  
\bean
	\E d_{\cF_0} (\bbP_n, P_0) \lesssim
	\left\{\begin{array}{ll}
	    n^{-\alpha/D} & ~~ \text{if $\alpha < D/2$}
	    \\
	    n^{-1/2} \log n & ~~ \text{if $\alpha = D/2$}
	    \\
	    n^{-1/2} & ~~ \text{if $\alpha > D/2$},
	\end{array}\right.
\eean
see \cite{schreuder2021bounding}, for example.
Similar bounds can be obtained for more general Besov IPMs.
Hence, when $\alpha/D < \beta_*/(2\beta_* + t_*) < 1/2$, the rate provided by the right-hand side of \eqref{eq:rate-ipm} is strictly faster than that of \eqref{eq:rate-ipm2}.

\section{Lower bound of the minimax risk} \label{sec:lower-bound}

In this section, we study a lower bound for the minimax optimal rate, particularly focusing on the case $d_{\rm eval} = W_1$. With $P_Z$ the uniform distribution on $[0,1]^d$, we investigate the minimax optimal rate for the distribution class $\cQ_0 = \{Q_\bg: \bg \in \cG_0\}$, where $\cG_0 = \cG_0(q, \bd, \bt, \bbeta, K)$. Our analysis is focused on the regime where $t_i \leq \min\{d_0, \ldots, d_{q+1}\}$ and $\beta_i \geq 1$ for all $i$. Beyond this regime, obtaining a lower bound using our proof technique becomes challenging. Note that $\widetilde \beta_i = \beta_i$ for all $i$ in this regime.

For given $\cG_0$ and $\sigma_0 \geq 0$, the minimax risk is defined as
\bean
	\mathfrak{M}(\cG_0, \sigma_0) = \inf_{\hat Q} \sup_{\bg_0 \in \cG_0} \E W_1(\hat Q, Q_0),
\eean
where the infimum ranges over all possible estimators.
Although the exact value of $\mathfrak{M}(\cG_0, \sigma_0)$ is rarely available in nonparametric problems, several techniques are known in the literature to obtain a lower bound of it.
We refer to \cite{tsybakov2008introduction} and \cite{wainwright2019high} for a comprehensive review.
We will utilize a general technique known as Fano's method to obtain a lower bound.

\begin{theorem} \label{thm:lower-bound}
Suppose that $d \leq D$, $\sigma_0 \geq 0$, $t_i \leq \min\{d_0, \ldots, d_{q+1}\}$, $\beta_i \geq 1$ for all $i$, $P_Z$ is the uniform distribution on $[0,1]^d$ and $\cG_0 = \cG_0(q, \bd, \bt, \bbeta, K)$.
If $K$ is large enough (depending on $\bbeta$ and $\bd$), there exists a constant $C > 0$ such that
\normalsize
\begin{align} \label{eq:lower-bound}
	\mathfrak{M}(\cG_0, \sigma_0) \geq C \max_{i \in \{0, \ldots, q\}} n^{ -\frac{\beta_i}{2\beta_i + t_i - 2}}.    
\end{align}
\end{theorem}

Note that the lower bound \eqref{eq:lower-bound} does not depend on $\sigma_0$.
With a direct application of Le Cam's method, one can easily show that $\mathfrak{M}(\cG_0, \sigma_0) \gtrsim \sigma_0 / \sqrt{n}$,  hence
\bean
	\mathfrak{M}(\cG_0, \sigma_0) \geq C \Big\{ \max_{i \in \{0, \ldots, q\}} n^{-\frac{\beta_i}{2\beta_i + t_i - 2}} + \frac{\sigma_0}{\sqrt{n}} \Big\}.
\eean
Since we are particularly interested in the small $\sigma_0$ regime (\ie\ nearly singular cases), our discussion below focuses on the case where $\sigma_0$ is small enough.

Firstly, note that the rate in the right hand side of \eqref{eq:lower-bound} can be strictly larger than $n^{-\beta_* / (2\beta_* + t_* - 2)}$, where $t_*$ and $\beta_*$ are defined in \eqref{eq:intrinsic-def}. If $(\beta_*, t_*) = (1, 2)$ and $(\beta_i, t_i) = (1.6, 3)$, for example, then $t_* / \beta_* > t_i / \beta_i$ and $\beta_i / (2\beta_i + t_i - 2) < \beta_* / (2\beta_* + t_* - 2)$. However, the rate in \eqref{eq:lower-bound} cannot be larger than $n^{-\beta_* / (2\beta_* + t_*)}$, which is the convergence rate (up to a logarithmic factor) in Theorem \ref{thm:rate-composition}. To provide a more convenient comparison of the upper and lower bounds of the convergence rate, we can express the bounds as follows
\bean
    &\text{Upper bound in Theorem \ref{thm:rate-composition}:}& \quad \max_{i\in \{0, \ldots, q\}} n^{-\frac{\beta_i}{2\beta_i + t_i}}
    \\
    &\text{Lower bound in Theorem \ref{thm:lower-bound}:}& \quad \max_{i \in \{0, \ldots, q\}} n^{-\frac{\beta_i}{2\beta_i + t_i - 2}}.
\eean
Therefore, the lower bound is only slightly smaller than the upper bound, indicating that the convergence rate of a GAN-based estimator is very close to the minimax optimal rate.

Regarding the difference between the upper and lower bounds, we conjecture that the lower bound is sharp in at least some special cases, and thus cannot be improved in general. In particular, when $q=0$ and $t_0=d=D$, we believe that the lower bound in Theorem \ref{thm:lower-bound} is sharp. This conjecture is based on the results presented in \cite{uppal2019nonparametric} and \cite{liang2021well}. 
They considered GAN for nonparametric density estimation, \ie\ $D=d$ in their framework.
For example, Theorem 4 in \cite{liang2021well} guarantees that, for $D=d\geq 2$ and $\sigma_0=0$, 
\be\label{eq:liang}
	\inf_{\hat Q} \sup_{Q_0 \in \widetilde\cQ_0} \E W_1(\hat Q, Q_0) \asymp n^{-\frac{\beta'+1}{2\beta'+d}},
\ee
where 
\bean
  \widetilde\cQ_0 = \Big\{Q: q \in \cH^{\beta'}_1([0,1]^d) \Big\}.
\eean
More precisely, they considered Sobolev classes rather than \Holder classes.
We also refer to \cite{niles2022minimax} for similar results but using different proof techniques.
Interestingly, when $D=d$, there is a close connection between the density model $\widetilde \cQ_0$ and the generative model $\cQ_0 = \{Q_\bg: \bg \in \cH_K^\beta([0,1]^d) \}$.
This connection is based on Caffarelli's regularity theory of optimal transport, often referred to as the Brenier map.
Roughly speaking, for a certain $\beta'$-\Holder density $q$, there exists a $(\beta'+1)$-\Holder function $\bg$ such that $Q = Q_\bg$. 
Therefore, a density model consisting of densities with this property can be viewed as a sub-model of the generative model with $\beta = \beta' +1$.
It is noteworthy that the rate \eqref{eq:liang} is the same as the right-hand side of \eqref{eq:lower-bound} with $q=0$, $d=D=t_0$, and $\beta_0=\beta=\beta'+1$.
This is why we conjecture that the lower bound \eqref{eq:lower-bound} cannot be improved in general.
It is important to note that this argument is a conjecture because Caffarelli's regularity theory requires the uniform convexity of the domain and co-domain of $\bg$, but $[0,1]^d$ is not uniformly convex.
For rigorous statements, counterexamples, and historical background on this field, we refer readers to Chapter 12 of \cite{villani2008optimal}.
Additionally, \cite{cordero2019regularity} provides some recent advancements in this area.

The minimax optimal rate in \cite{tang2023minimax} is consistent with the lower bound \eqref{eq:lower-bound} in some special cases with $q=0$. However, the proof techniques presented in the existing literature are not directly applicable to the structured distribution estimation problem considered in our paper. The techniques employed in \cite{uppal2019nonparametric} and \cite{liang2021well} for both upper and lower bounds rely on wavelet thresholding. It is unclear how to extend these techniques to our case with $q > 0$ and a singular $Q_0$. Additionally, \cite{tang2023minimax} also relies on wavelet thresholding for estimating a distribution on a manifold, which is similar to the techniques used in \cite{uppal2019nonparametric} and \cite{liang2021well}, making them not fully applicable to our specific estimation problem. Although \cite{tang2023minimax} incorporates an additional step of estimating the charts of the manifold, these techniques do not directly address our estimation problem.

\section{Conclusion} \label{sec:conclusion}

Under a structural assumption on the generator, we investigated a convergence rate of a GAN-based estimator and a lower bound of the minimax optimal rate.
Notably, the rate is faster than that obtained by likelihood approaches.
In practice, however, the computation of GAN incorporates a challenging minimax optimization problem and our understanding of it remains largely unexplored.
For example, it is unclear where a practical estimator constructed via a stochastic gradient algorithm converges to  \citep{mescheder2018training, hsieh2021limits}.
The discriminator constructed in the proof of Theorem \ref{thm:rate-composition} is even further away from the one used in practice.
Our theory only guarantees that there exists a discriminator class $\cF$ which yields an estimator whose convergence rate is close to the minimax optimal rate.
Regardless, our theory plays an important role in further advancing GAN theory.

We conclude the paper with some possible directions for future work. One of the most important tasks is to reduce the current gap between the upper and lower bounds of the convergence rate. As discussed in Section \ref{sec:lower-bound}, it would be crucial to construct an estimator that achieves the lower bound in Theorem \ref{thm:lower-bound}. 
After an early version of this paper was drafted on arXiv, \cite{stephanovitch2023wasserstein} studied a similar problem, particularly focusing on the special case where $q=0$ and $t_0=d$. They obtained a minimax optimal estimator, but its construction relies on wavelet features rather than DNNs. Furthermore, their proof techniques cannot be extended to more general cases where $q > 0$. Techniques from the literature concerning the estimation of optimal transport maps might be employed to address this problem, as explored in works such as 
\cite{deb2021rates},
\cite{hutter2021minimax},
\cite{divol2022optimal},
\cite{manole2021plugin} and
\cite{pooladian2021entropic}. The problem of estimating optimal transport maps appears to be closely related to our set-up, and the rate \eqref{thm:lower-bound} can be found in this literature. Investigating whether a GAN-based estimator can achieve the minimax rate is another important research problem. In particular, it would be valuable to explore whether the discriminator and generator classes modeled by deep neural networks can attain the minimax rate when $d_{\rm eval} = W_1$. Finally, based on the approximation property of the convolutional neural networks (CNN) architectures \citep{kohler2020rate, yarotsky2021universal}, studying the benefit of CNN-based GAN would be an intriguing problem.

\section*{Acknowledgments}
The author would like to thank Lizhen Lin for her valuable comments and discussions on an earlier version of this paper. This work was supported by Samsung Science and Technology Foundation under Project Number SSTF-BA2101-03.

\bibliography{sn-bibliography}

\begin{thebibliography}{}
\renewcommand{\doi}[1]{\url{https://doi.org/#1}}
\bibcommenthead

\bibitem [\protect \citeauthoryear {%
Aamari%
\ \BBA {} Levrard%
}{%
Aamari%
\ \BBA {} Levrard%
}{%
{\protect \APACyear {2018}}%
}]{%
aamari2018stability}
\APACinsertmetastar {%
aamari2018stability}%
\begin{APACrefauthors}%
Aamari, E.%
\BCBT {}\ \BBA {} Levrard, C.%
\end{APACrefauthors}%
\unskip\
\newblock
\APACrefYearMonthDay{2018}{}{}.
\newblock
{\BBOQ}\APACrefatitle {{Stability and minimax optimality of tangential Delaunay complexes for manifold reconstruction}} {{Stability and minimax optimality of tangential Delaunay complexes for manifold reconstruction}}.{\BBCQ}
\newblock
\APACjournalVolNumPages{Discrete Comput. Geom.}{59}{4}{923--971,}
\newblock

\newblock

\PrintBackRefs{\CurrentBib}

\bibitem [\protect \citeauthoryear {%
Aamari%
\ \BBA {} Levrard%
}{%
Aamari%
\ \BBA {} Levrard%
}{%
{\protect \APACyear {2019}}%
}]{%
aamari2019nonasymptotic}
\APACinsertmetastar {%
aamari2019nonasymptotic}%
\begin{APACrefauthors}%
Aamari, E.%
\BCBT {}\ \BBA {} Levrard, C.%
\end{APACrefauthors}%
\unskip\
\newblock
\APACrefYearMonthDay{2019}{}{}.
\newblock
{\BBOQ}\APACrefatitle {Nonasymptotic rates for manifold, tangent space and curvature estimation} {Nonasymptotic rates for manifold, tangent space and curvature estimation}.{\BBCQ}
\newblock
\APACjournalVolNumPages{Ann. Statist.}{47}{1}{177--204,}
\newblock

\newblock

\PrintBackRefs{\CurrentBib}

\bibitem [\protect \citeauthoryear {%
Arjovsky%
, Chintala%
\BCBL {}\ \BBA {} Bottou%
}{%
Arjovsky%
\ \protect \BOthers {.}}{%
{\protect \APACyear {2017}}%
}]{%
arjovsky2017wasserstein}
\APACinsertmetastar {%
arjovsky2017wasserstein}%
\begin{APACrefauthors}%
Arjovsky, M.%
, Chintala, S.%
\BCBL {} Bottou, L.%
\end{APACrefauthors}%
\unskip\
\newblock
\APACrefYearMonthDay{2017}{}{}.
\newblock
{\BBOQ}\APACrefatitle {Wasserstein generative adversarial networks} {Wasserstein generative adversarial networks}.{\BBCQ}
\newblock
 \APACrefbtitle {{Proc. ICML}} {{Proc. ICML}}\ (\BPGS\ 214--223).
\PrintBackRefs{\CurrentBib}

\bibitem [\protect \citeauthoryear {%
Arora%
, Ge%
, Liang%
, Ma%
\BCBL {}\ \BBA {} Zhang%
}{%
Arora%
\ \protect \BOthers {.}}{%
{\protect \APACyear {2017}}%
}]{%
arora2017generalization}
\APACinsertmetastar {%
arora2017generalization}%
\begin{APACrefauthors}%
Arora, S.%
, Ge, R.%
, Liang, Y.%
, Ma, T.%
\BCBL {} Zhang, Y.%
\end{APACrefauthors}%
\unskip\
\newblock
\APACrefYearMonthDay{2017}{}{}.
\newblock
{\BBOQ}\APACrefatitle {{Generalization and equilibrium in generative adversarial nets (GANs)}} {{Generalization and equilibrium in generative adversarial nets (GANs)}}.{\BBCQ}
\newblock
 \APACrefbtitle {{Proc. ICML}} {{Proc. ICML}}\ (\BPGS\ 224--232).
\PrintBackRefs{\CurrentBib}

\bibitem [\protect \citeauthoryear {%
Bai%
, Ma%
\BCBL {}\ \BBA {} Risteski%
}{%
Bai%
\ \protect \BOthers {.}}{%
{\protect \APACyear {2019}}%
}]{%
bai2019approximability}
\APACinsertmetastar {%
bai2019approximability}%
\begin{APACrefauthors}%
Bai, Y.%
, Ma, T.%
\BCBL {} Risteski, A.%
\end{APACrefauthors}%
\unskip\
\newblock
\APACrefYearMonthDay{2019}{}{}.
\newblock
{\BBOQ}\APACrefatitle {{Approximability of discriminators implies diversity in GANs}} {{Approximability of discriminators implies diversity in GANs}}.{\BBCQ}
\newblock
 \APACrefbtitle {{Proc. ICLR}} {{Proc. ICLR}}\ (\BPGS\ 1--10).
\PrintBackRefs{\CurrentBib}

\bibitem [\protect \citeauthoryear {%
Bauer%
\ \BBA {} Kohler%
}{%
Bauer%
\ \BBA {} Kohler%
}{%
{\protect \APACyear {2019}}%
}]{%
bauer2019deep}
\APACinsertmetastar {%
bauer2019deep}%
\begin{APACrefauthors}%
Bauer, B.%
\BCBT {}\ \BBA {} Kohler, M.%
\end{APACrefauthors}%
\unskip\
\newblock
\APACrefYearMonthDay{2019}{}{}.
\newblock
{\BBOQ}\APACrefatitle {On deep learning as a remedy for the curse of dimensionality in nonparametric regression} {On deep learning as a remedy for the curse of dimensionality in nonparametric regression}.{\BBCQ}
\newblock
\APACjournalVolNumPages{Ann. Statist.}{47}{4}{2261--2285,}
\newblock

\newblock

\PrintBackRefs{\CurrentBib}

\bibitem [\protect \citeauthoryear {%
Belomestny%
, Moulines%
, Naumov%
, Puchkin%
\BCBL {}\ \BBA {} Samsonov%
}{%
Belomestny%
\ \protect \BOthers {.}}{%
{\protect \APACyear {2021}}%
}]{%
belomestny2021rates}
\APACinsertmetastar {%
belomestny2021rates}%
\begin{APACrefauthors}%
Belomestny, D.%
, Moulines, E.%
, Naumov, A.%
, Puchkin, N.%
\BCBL {} Samsonov, S.%
\end{APACrefauthors}%
\unskip\
\newblock
\APACrefYearMonthDay{2021}{}{}.
\newblock
{\BBOQ}\APACrefatitle {{Rates of convergence for density estimation with GANs}} {{Rates of convergence for density estimation with GANs}}.{\BBCQ}
\newblock
\APACjournalVolNumPages{ArXiv:2102.00199}{}{}{,}
\newblock

\newblock

\PrintBackRefs{\CurrentBib}

\bibitem [\protect \citeauthoryear {%
Berenfeld%
\ \BBA {} Hoffmann%
}{%
Berenfeld%
\ \BBA {} Hoffmann%
}{%
{\protect \APACyear {2021}}%
}]{%
berenfeld2021density}
\APACinsertmetastar {%
berenfeld2021density}%
\begin{APACrefauthors}%
Berenfeld, C.%
\BCBT {}\ \BBA {} Hoffmann, M.%
\end{APACrefauthors}%
\unskip\
\newblock
\APACrefYearMonthDay{2021}{}{}.
\newblock
{\BBOQ}\APACrefatitle {Density estimation on an unknown submanifold} {Density estimation on an unknown submanifold}.{\BBCQ}
\newblock
\APACjournalVolNumPages{Electron. J. Stat.}{15}{}{2179--2223,}
\newblock

\newblock

\PrintBackRefs{\CurrentBib}

\bibitem [\protect \citeauthoryear {%
Berenfeld%
, Rosa%
\BCBL {}\ \BBA {} Rousseau%
}{%
Berenfeld%
\ \protect \BOthers {.}}{%
{\protect \APACyear {2022}}%
}]{%
berenfeld2022estimating}
\APACinsertmetastar {%
berenfeld2022estimating}%
\begin{APACrefauthors}%
Berenfeld, C.%
, Rosa, P.%
\BCBL {} Rousseau, J.%
\end{APACrefauthors}%
\unskip\
\newblock
\APACrefYearMonthDay{2022}{}{}.
\newblock
{\BBOQ}\APACrefatitle {{Estimating a density near an unknown manifold: a Bayesian nonparametric approach}} {{Estimating a density near an unknown manifold: a Bayesian nonparametric approach}}.{\BBCQ}
\newblock
\APACjournalVolNumPages{ArXiv:2205.15717}{}{}{,}
\newblock

\newblock

\PrintBackRefs{\CurrentBib}

\bibitem [\protect \citeauthoryear {%
Biau%
, Sangnier%
\BCBL {}\ \BBA {} Tanielian%
}{%
Biau%
\ \protect \BOthers {.}}{%
{\protect \APACyear {2021}}%
}]{%
biau2021some}
\APACinsertmetastar {%
biau2021some}%
\begin{APACrefauthors}%
Biau, G.%
, Sangnier, M.%
\BCBL {} Tanielian, U.%
\end{APACrefauthors}%
\unskip\
\newblock
\APACrefYearMonthDay{2021}{}{}.
\newblock
{\BBOQ}\APACrefatitle {{Some theoretical insights into Wasserstein GANs}} {{Some theoretical insights into Wasserstein GANs}}.{\BBCQ}
\newblock
\APACjournalVolNumPages{J. Mach. Learn. Res.}{22}{1}{5287--5331,}
\newblock

\newblock

\PrintBackRefs{\CurrentBib}

\bibitem [\protect \citeauthoryear {%
B{\"u}hlmann%
\ \BBA {} van~de Geer%
}{%
B{\"u}hlmann%
\ \BBA {} van~de Geer%
}{%
{\protect \APACyear {2011}}%
}]{%
buhlmann2011statistics}
\APACinsertmetastar {%
buhlmann2011statistics}%
\begin{APACrefauthors}%
B{\"u}hlmann, P.%
\BCBT {}\ \BBA {} van~de Geer, S.%
\end{APACrefauthors}%
\unskip\
\newblock
\APACrefYear{2011}.
\newblock
\APACrefbtitle {{Statistics for High-Dimensional Data: Methods, Theory and Applications}} {{Statistics for High-Dimensional Data: Methods, Theory and Applications}}.
\newblock
\APACaddressPublisher{}{Springer}.
\PrintBackRefs{\CurrentBib}

\bibitem [\protect \citeauthoryear {%
Caffarelli%
}{%
Caffarelli%
}{%
{\protect \APACyear {1990}}%
}]{%
caffarelli1990interior}
\APACinsertmetastar {%
caffarelli1990interior}%
\begin{APACrefauthors}%
Caffarelli, L.A.%
\end{APACrefauthors}%
\unskip\
\newblock
\APACrefYearMonthDay{1990}{}{}.
\newblock
{\BBOQ}\APACrefatitle {{Interior $W^{2,p}$ estimates for solutions of the Monge--Amp\`{e}re equation}} {{Interior $W^{2,p}$ estimates for solutions of the Monge--Amp\`{e}re equation}}.{\BBCQ}
\newblock
\APACjournalVolNumPages{Ann. of Math.}{131}{1}{135--150,}
\newblock

\newblock

\PrintBackRefs{\CurrentBib}

\bibitem [\protect \citeauthoryear {%
Chae%
, Kim%
, Kim%
\BCBL {}\ \BBA {} Lin%
}{%
Chae%
\ \protect \BOthers {.}}{%
{\protect \APACyear {2023}}%
}]{%
chae2023likelihood}
\APACinsertmetastar {%
chae2023likelihood}%
\begin{APACrefauthors}%
Chae, M.%
, Kim, D.%
, Kim, Y.%
\BCBL {} Lin, L.%
\end{APACrefauthors}%
\unskip\
\newblock
\APACrefYearMonthDay{2023}{}{}.
\newblock
{\BBOQ}\APACrefatitle {A likelihood approach to nonparametric estimation of a singular distribution using deep generative models} {A likelihood approach to nonparametric estimation of a singular distribution using deep generative models}.{\BBCQ}
\newblock
\APACjournalVolNumPages{J. Mach. Learn. Res.}{24}{77}{1--42,}
\newblock

\newblock

\PrintBackRefs{\CurrentBib}

\bibitem [\protect \citeauthoryear {%
Chae%
\ \BBA {} Walker%
}{%
Chae%
\ \BBA {} Walker%
}{%
{\protect \APACyear {2019}}%
}]{%
chae2019bayesianConsistency}
\APACinsertmetastar {%
chae2019bayesianConsistency}%
\begin{APACrefauthors}%
Chae, M.%
\BCBT {}\ \BBA {} Walker, S.G.%
\end{APACrefauthors}%
\unskip\
\newblock
\APACrefYearMonthDay{2019}{}{}.
\newblock
{\BBOQ}\APACrefatitle {{Bayesian consistency for a nonparametric stationary Markov model}} {{Bayesian consistency for a nonparametric stationary Markov model}}.{\BBCQ}
\newblock
\APACjournalVolNumPages{Bernoulli}{25}{2}{877--901,}
\newblock

\newblock

\PrintBackRefs{\CurrentBib}

\bibitem [\protect \citeauthoryear {%
Chen%
, Liao%
, Zha%
\BCBL {}\ \BBA {} Zhao%
}{%
Chen%
\ \protect \BOthers {.}}{%
{\protect \APACyear {2020}}%
}]{%
chen2020statistical}
\APACinsertmetastar {%
chen2020statistical}%
\begin{APACrefauthors}%
Chen, M.%
, Liao, W.%
, Zha, H.%
\BCBL {} Zhao, T.%
\end{APACrefauthors}%
\unskip\
\newblock
\APACrefYearMonthDay{2020}{}{}.
\newblock
{\BBOQ}\APACrefatitle {Statistical guarantees of generative adversarial networks for distribution estimation} {Statistical guarantees of generative adversarial networks for distribution estimation}.{\BBCQ}
\newblock
\APACjournalVolNumPages{ArXiv:2002.03938}{}{}{,}
\newblock

\newblock

\PrintBackRefs{\CurrentBib}

\bibitem [\protect \citeauthoryear {%
Cordero-Erausquin%
\ \BBA {} Figalli%
}{%
Cordero-Erausquin%
\ \BBA {} Figalli%
}{%
{\protect \APACyear {2019}}%
}]{%
cordero2019regularity}
\APACinsertmetastar {%
cordero2019regularity}%
\begin{APACrefauthors}%
Cordero-Erausquin, D.%
\BCBT {}\ \BBA {} Figalli, A.%
\end{APACrefauthors}%
\unskip\
\newblock
\APACrefYearMonthDay{2019}{}{}.
\newblock
{\BBOQ}\APACrefatitle {Regularity of monotone transport maps between unbounded domains} {Regularity of monotone transport maps between unbounded domains}.{\BBCQ}
\newblock
\APACjournalVolNumPages{Discrete Contin. Dyn. Syst.}{39}{12}{7101-7112,}
\newblock

\newblock

\PrintBackRefs{\CurrentBib}

\bibitem [\protect \citeauthoryear {%
Deb%
, Ghosal%
\BCBL {}\ \BBA {} Sen%
}{%
Deb%
\ \protect \BOthers {.}}{%
{\protect \APACyear {2021}}%
}]{%
deb2021rates}
\APACinsertmetastar {%
deb2021rates}%
\begin{APACrefauthors}%
Deb, N.%
, Ghosal, P.%
\BCBL {} Sen, B.%
\end{APACrefauthors}%
\unskip\
\newblock
\APACrefYearMonthDay{2021}{}{}.
\newblock
{\BBOQ}\APACrefatitle {Rates of estimation of optimal transport maps using plug-in estimators via barycentric projections} {Rates of estimation of optimal transport maps using plug-in estimators via barycentric projections}.{\BBCQ}
\newblock
 \APACrefbtitle {{Proc. NeurIPS}} {{Proc. NeurIPS}}\ (\BVOL~34, \BPGS\ 29736--29753).
\PrintBackRefs{\CurrentBib}

\bibitem [\protect \citeauthoryear {%
Divol%
}{%
Divol%
}{%
{\protect \APACyear {2020}}%
}]{%
divol2020minimax}
\APACinsertmetastar {%
divol2020minimax}%
\begin{APACrefauthors}%
Divol, V.%
\end{APACrefauthors}%
\unskip\
\newblock
\APACrefYearMonthDay{2020}{}{}.
\newblock
{\BBOQ}\APACrefatitle {Minimax adaptive estimation in manifold inference} {Minimax adaptive estimation in manifold inference}.{\BBCQ}
\newblock
\APACjournalVolNumPages{ArXiv:2001.04896}{}{}{,}
\newblock

\newblock

\PrintBackRefs{\CurrentBib}

\bibitem [\protect \citeauthoryear {%
Divol%
}{%
Divol%
}{%
{\protect \APACyear {2021}}%
}]{%
divol2021minimax}
\APACinsertmetastar {%
divol2021minimax}%
\begin{APACrefauthors}%
Divol, V.%
\end{APACrefauthors}%
\unskip\
\newblock
\APACrefYearMonthDay{2021}{}{}.
\newblock
{\BBOQ}\APACrefatitle {Minimax adaptive estimation in manifold inference} {Minimax adaptive estimation in manifold inference}.{\BBCQ}
\newblock
\APACjournalVolNumPages{Electron. J. Stat.}{15}{2}{5888--5932,}
\newblock

\newblock

\PrintBackRefs{\CurrentBib}

\bibitem [\protect \citeauthoryear {%
Divol%
}{%
Divol%
}{%
{\protect \APACyear {2022}}%
}]{%
divol2022measure}
\APACinsertmetastar {%
divol2022measure}%
\begin{APACrefauthors}%
Divol, V.%
\end{APACrefauthors}%
\unskip\
\newblock
\APACrefYearMonthDay{2022}{}{}.
\newblock
{\BBOQ}\APACrefatitle {Measure estimation on manifolds: an optimal transport approach} {Measure estimation on manifolds: an optimal transport approach}.{\BBCQ}
\newblock
\APACjournalVolNumPages{Probab. Theory Related Fields}{183}{1}{581--647,}
\newblock

\newblock

\PrintBackRefs{\CurrentBib}

\bibitem [\protect \citeauthoryear {%
Divol%
, Niles-Weed%
\BCBL {}\ \BBA {} Pooladian%
}{%
Divol%
\ \protect \BOthers {.}}{%
{\protect \APACyear {2022}}%
}]{%
divol2022optimal}
\APACinsertmetastar {%
divol2022optimal}%
\begin{APACrefauthors}%
Divol, V.%
, Niles-Weed, J.%
\BCBL {} Pooladian, A\BHBI A.%
\end{APACrefauthors}%
\unskip\
\newblock
\APACrefYearMonthDay{2022}{}{}.
\newblock
{\BBOQ}\APACrefatitle {Optimal transport map estimation in general function spaces} {Optimal transport map estimation in general function spaces}.{\BBCQ}
\newblock
\APACjournalVolNumPages{ArXiv:2212.03722}{}{}{,}
\newblock

\newblock

\PrintBackRefs{\CurrentBib}

\bibitem [\protect \citeauthoryear {%
Fan%
}{%
Fan%
}{%
{\protect \APACyear {1991}}%
}]{%
fan1991optimal}
\APACinsertmetastar {%
fan1991optimal}%
\begin{APACrefauthors}%
Fan, J.%
\end{APACrefauthors}%
\unskip\
\newblock
\APACrefYearMonthDay{1991}{}{}.
\newblock
{\BBOQ}\APACrefatitle {On the optimal rates of convergence for nonparametric deconvolution problems} {On the optimal rates of convergence for nonparametric deconvolution problems}.{\BBCQ}
\newblock
\APACjournalVolNumPages{Ann. Statist.}{19}{3}{1257--1272,}
\newblock

\newblock

\PrintBackRefs{\CurrentBib}

\bibitem [\protect \citeauthoryear {%
Fournier%
\ \BBA {} Guillin%
}{%
Fournier%
\ \BBA {} Guillin%
}{%
{\protect \APACyear {2015}}%
}]{%
fournier2015rate}
\APACinsertmetastar {%
fournier2015rate}%
\begin{APACrefauthors}%
Fournier, N.%
\BCBT {}\ \BBA {} Guillin, A.%
\end{APACrefauthors}%
\unskip\
\newblock
\APACrefYearMonthDay{2015}{}{}.
\newblock
{\BBOQ}\APACrefatitle {{On the rate of convergence in Wasserstein distance of the empirical measure}} {{On the rate of convergence in Wasserstein distance of the empirical measure}}.{\BBCQ}
\newblock
\APACjournalVolNumPages{Probab. Theory Related Fields}{162}{3-4}{707--738,}
\newblock

\newblock

\PrintBackRefs{\CurrentBib}

\bibitem [\protect \citeauthoryear {%
Genovese%
, Perone-Pacifico%
, Verdinelli%
\BCBL {}\ \BBA {} Wasserman%
}{%
Genovese%
\ \protect \BOthers {.}}{%
{\protect \APACyear {2012}}%
{\protect \APACexlab {{\protect \BCnt {1}}}}}]{%
genovese2012manifold}
\APACinsertmetastar {%
genovese2012manifold}%
\begin{APACrefauthors}%
Genovese, C.R.%
, Perone-Pacifico, M.%
, Verdinelli, I.%
\BCBL {} Wasserman, L.%
\end{APACrefauthors}%
\unskip\
\newblock
\APACrefYearMonthDay{2012{\protect \BCnt {1}}}{}{}.
\newblock
{\BBOQ}\APACrefatitle {{Manifold estimation and singular deconvolution under Hausdorff loss}} {{Manifold estimation and singular deconvolution under Hausdorff loss}}.{\BBCQ}
\newblock
\APACjournalVolNumPages{Ann. Statist.}{40}{2}{941--963,}
\newblock

\newblock

\PrintBackRefs{\CurrentBib}

\bibitem [\protect \citeauthoryear {%
Genovese%
, Perone-Pacifico%
, Verdinelli%
\BCBL {}\ \BBA {} Wasserman%
}{%
Genovese%
\ \protect \BOthers {.}}{%
{\protect \APACyear {2012}}%
{\protect \APACexlab {{\protect \BCnt {2}}}}}]{%
genovese2012minimax}
\APACinsertmetastar {%
genovese2012minimax}%
\begin{APACrefauthors}%
Genovese, C.R.%
, Perone-Pacifico, M.%
, Verdinelli, I.%
\BCBL {} Wasserman, L.%
\end{APACrefauthors}%
\unskip\
\newblock
\APACrefYearMonthDay{2012{\protect \BCnt {2}}}{}{}.
\newblock
{\BBOQ}\APACrefatitle {Minimax manifold estimation} {Minimax manifold estimation}.{\BBCQ}
\newblock
\APACjournalVolNumPages{J. Mach. Learn. Res.}{13}{1}{1263--1291,}
\newblock

\newblock

\PrintBackRefs{\CurrentBib}

\bibitem [\protect \citeauthoryear {%
Ghosal%
\ \BBA {} van~der Vaart%
}{%
Ghosal%
\ \BBA {} van~der Vaart%
}{%
{\protect \APACyear {2017}}%
}]{%
ghosal2017fundamentals}
\APACinsertmetastar {%
ghosal2017fundamentals}%
\begin{APACrefauthors}%
Ghosal, S.%
\BCBT {}\ \BBA {} van~der Vaart, A.%
\end{APACrefauthors}%
\unskip\
\newblock
\APACrefYear{2017}.
\newblock
\APACrefbtitle {{Fundamentals of Nonparametric Bayesian Inference}} {{Fundamentals of Nonparametric Bayesian Inference}}.
\newblock
\APACaddressPublisher{}{Cambridge University Press}.
\PrintBackRefs{\CurrentBib}

\bibitem [\protect \citeauthoryear {%
Gin{\'e}%
\ \BBA {} Nickl%
}{%
Gin{\'e}%
\ \BBA {} Nickl%
}{%
{\protect \APACyear {2016}}%
}]{%
gine2016mathematical}
\APACinsertmetastar {%
gine2016mathematical}%
\begin{APACrefauthors}%
Gin{\'e}, E.%
\BCBT {}\ \BBA {} Nickl, R.%
\end{APACrefauthors}%
\unskip\
\newblock
\APACrefYear{2016}.
\newblock
\APACrefbtitle {Mathematical Foundations of Infinite-Dimensional Statistical Models} {Mathematical foundations of infinite-dimensional statistical models}.
\newblock
\APACaddressPublisher{}{Cambridge University Press}.
\PrintBackRefs{\CurrentBib}

\bibitem [\protect \citeauthoryear {%
Givens%
\ \BBA {} Shortt%
}{%
Givens%
\ \BBA {} Shortt%
}{%
{\protect \APACyear {1984}}%
}]{%
givens1984class}
\APACinsertmetastar {%
givens1984class}%
\begin{APACrefauthors}%
Givens, C.R.%
\BCBT {}\ \BBA {} Shortt, R.M.%
\end{APACrefauthors}%
\unskip\
\newblock
\APACrefYearMonthDay{1984}{}{}.
\newblock
{\BBOQ}\APACrefatitle {{A class of Wasserstein metrics for probability distributions}} {{A class of Wasserstein metrics for probability distributions}}.{\BBCQ}
\newblock
\APACjournalVolNumPages{The Michigan Mathematical Journal}{31}{2}{231--240,}
\newblock

\newblock

\PrintBackRefs{\CurrentBib}

\bibitem [\protect \citeauthoryear {%
Glorot%
, Bordes%
\BCBL {}\ \BBA {} Bengio%
}{%
Glorot%
\ \protect \BOthers {.}}{%
{\protect \APACyear {2011}}%
}]{%
glorot2011deep}
\APACinsertmetastar {%
glorot2011deep}%
\begin{APACrefauthors}%
Glorot, X.%
, Bordes, A.%
\BCBL {} Bengio, Y.%
\end{APACrefauthors}%
\unskip\
\newblock
\APACrefYearMonthDay{2011}{}{}.
\newblock
{\BBOQ}\APACrefatitle {Deep sparse rectifier neural networks} {Deep sparse rectifier neural networks}.{\BBCQ}
\newblock
 \APACrefbtitle {{Proc. AISTATS}} {{Proc. AISTATS}}\ (\BPGS\ 315--323).
\PrintBackRefs{\CurrentBib}

\bibitem [\protect \citeauthoryear {%
Goodfellow%
\ \protect \BOthers {.}}{%
Goodfellow%
\ \protect \BOthers {.}}{%
{\protect \APACyear {2014}}%
}]{%
goodfellow2014generative}
\APACinsertmetastar {%
goodfellow2014generative}%
\begin{APACrefauthors}%
Goodfellow, I.%
, Pouget-Abadie, J.%
, Mirza, M.%
, Xu, B.%
, Warde-Farley, D.%
, Ozair, S.%
\BDBL {}Bengio, Y.%
\end{APACrefauthors}%
\unskip\
\newblock
\APACrefYearMonthDay{2014}{}{}.
\newblock
{\BBOQ}\APACrefatitle {Generative adversarial nets} {Generative adversarial nets}.{\BBCQ}
\newblock
 \APACrefbtitle {{Proc. NIPS}} {{Proc. NIPS}}\ (\BPGS\ 2672--2680).
\PrintBackRefs{\CurrentBib}

\bibitem [\protect \citeauthoryear {%
Hastie%
, Tibshirani%
\BCBL {}\ \BBA {} Friedman%
}{%
Hastie%
\ \protect \BOthers {.}}{%
{\protect \APACyear {2009}}%
}]{%
hastie2009elements}
\APACinsertmetastar {%
hastie2009elements}%
\begin{APACrefauthors}%
Hastie, T.%
, Tibshirani, R.%
\BCBL {} Friedman, J.H.%
\end{APACrefauthors}%
\unskip\
\newblock
\APACrefYear{2009}.
\newblock
\APACrefbtitle {{The Elements of Statistical Learning: Data Mining, Inference, and Prediction}} {{The Elements of Statistical Learning: Data Mining, Inference, and Prediction}}.
\newblock
\APACaddressPublisher{}{Springer, New York}.
\PrintBackRefs{\CurrentBib}

\bibitem [\protect \citeauthoryear {%
Horowitz%
\ \BBA {} Mammen%
}{%
Horowitz%
\ \BBA {} Mammen%
}{%
{\protect \APACyear {2007}}%
}]{%
horowitz2007rate}
\APACinsertmetastar {%
horowitz2007rate}%
\begin{APACrefauthors}%
Horowitz, J.L.%
\BCBT {}\ \BBA {} Mammen, E.%
\end{APACrefauthors}%
\unskip\
\newblock
\APACrefYearMonthDay{2007}{}{}.
\newblock
{\BBOQ}\APACrefatitle {Rate-optimal estimation for a general class of nonparametric regression models with unknown link functions} {Rate-optimal estimation for a general class of nonparametric regression models with unknown link functions}.{\BBCQ}
\newblock
\APACjournalVolNumPages{Ann. Statist.}{35}{6}{2589--2619,}
\newblock

\newblock

\PrintBackRefs{\CurrentBib}

\bibitem [\protect \citeauthoryear {%
Hsieh%
, Mertikopoulos%
\BCBL {}\ \BBA {} Cevher%
}{%
Hsieh%
\ \protect \BOthers {.}}{%
{\protect \APACyear {2021}}%
}]{%
hsieh2021limits}
\APACinsertmetastar {%
hsieh2021limits}%
\begin{APACrefauthors}%
Hsieh, Y\BHBI P.%
, Mertikopoulos, P.%
\BCBL {} Cevher, V.%
\end{APACrefauthors}%
\unskip\
\newblock
\APACrefYearMonthDay{2021}{}{}.
\newblock
{\BBOQ}\APACrefatitle {The limits of min-max optimization algorithms: Convergence to spurious non-critical sets} {The limits of min-max optimization algorithms: Convergence to spurious non-critical sets}.{\BBCQ}
\newblock
 \APACrefbtitle {{Proc. ICML}} {{Proc. ICML}}\ (\BPGS\ 4337--4348).
\PrintBackRefs{\CurrentBib}

\bibitem [\protect \citeauthoryear {%
H{\"u}tter%
\ \BBA {} Rigollet%
}{%
H{\"u}tter%
\ \BBA {} Rigollet%
}{%
{\protect \APACyear {2021}}%
}]{%
hutter2021minimax}
\APACinsertmetastar {%
hutter2021minimax}%
\begin{APACrefauthors}%
H{\"u}tter, J\BHBI C.%
\BCBT {}\ \BBA {} Rigollet, P.%
\end{APACrefauthors}%
\unskip\
\newblock
\APACrefYearMonthDay{2021}{}{}.
\newblock
{\BBOQ}\APACrefatitle {Minimax estimation of smooth optimal transport maps} {Minimax estimation of smooth optimal transport maps}.{\BBCQ}
\newblock
\APACjournalVolNumPages{Ann. Statist.}{49}{}{1166--1194,}
\newblock

\newblock

\PrintBackRefs{\CurrentBib}

\bibitem [\protect \citeauthoryear {%
Jiao%
, Shen%
, Lin%
\BCBL {}\ \BBA {} Huang%
}{%
Jiao%
\ \protect \BOthers {.}}{%
{\protect \APACyear {2023}}%
}]{%
jiao2023deep}
\APACinsertmetastar {%
jiao2023deep}%
\begin{APACrefauthors}%
Jiao, Y.%
, Shen, G.%
, Lin, Y.%
\BCBL {} Huang, J.%
\end{APACrefauthors}%
\unskip\
\newblock
\APACrefYearMonthDay{2023}{}{}.
\newblock
{\BBOQ}\APACrefatitle {{Deep nonparametric regression on approximate manifolds: Nonasymptotic error bounds with polynomial prefactors}} {{Deep nonparametric regression on approximate manifolds: Nonasymptotic error bounds with polynomial prefactors}}.{\BBCQ}
\newblock
\APACjournalVolNumPages{Ann. Statist.}{}{2}{691 -- 716,}
\newblock

\newblock

\PrintBackRefs{\CurrentBib}

\bibitem [\protect \citeauthoryear {%
Juditsky%
, Lepski%
\BCBL {}\ \BBA {} Tsybakov%
}{%
Juditsky%
\ \protect \BOthers {.}}{%
{\protect \APACyear {2009}}%
}]{%
juditsky2009nonparametric}
\APACinsertmetastar {%
juditsky2009nonparametric}%
\begin{APACrefauthors}%
Juditsky, A.B.%
, Lepski, O.V.%
\BCBL {} Tsybakov, A.B.%
\end{APACrefauthors}%
\unskip\
\newblock
\APACrefYearMonthDay{2009}{}{}.
\newblock
{\BBOQ}\APACrefatitle {Nonparametric estimation of composite functions} {Nonparametric estimation of composite functions}.{\BBCQ}
\newblock
\APACjournalVolNumPages{Ann. Statist.}{37}{3}{1360--1404,}
\newblock

\newblock

\PrintBackRefs{\CurrentBib}

\bibitem [\protect \citeauthoryear {%
Karras%
, Aila%
, Laine%
\BCBL {}\ \BBA {} Lehtinen%
}{%
Karras%
\ \protect \BOthers {.}}{%
{\protect \APACyear {2018}}%
}]{%
karras2018progressive}
\APACinsertmetastar {%
karras2018progressive}%
\begin{APACrefauthors}%
Karras, T.%
, Aila, T.%
, Laine, S.%
\BCBL {} Lehtinen, J.%
\end{APACrefauthors}%
\unskip\
\newblock
\APACrefYearMonthDay{2018}{}{}.
\newblock
{\BBOQ}\APACrefatitle {{Progressive growing of GANs for improved quality, stability, and variation}} {{Progressive growing of GANs for improved quality, stability, and variation}}.{\BBCQ}
\newblock
 \APACrefbtitle {{Proc. ICLR}} {{Proc. ICLR}}\ (\BPGS\ 1--26).
\PrintBackRefs{\CurrentBib}

\bibitem [\protect \citeauthoryear {%
Karras%
, Laine%
\BCBL {}\ \BBA {} Aila%
}{%
Karras%
\ \protect \BOthers {.}}{%
{\protect \APACyear {2019}}%
}]{%
karras2019style}
\APACinsertmetastar {%
karras2019style}%
\begin{APACrefauthors}%
Karras, T.%
, Laine, S.%
\BCBL {} Aila, T.%
\end{APACrefauthors}%
\unskip\
\newblock
\APACrefYearMonthDay{2019}{}{}.
\newblock
{\BBOQ}\APACrefatitle {A style-based generator architecture for generative adversarial networks} {A style-based generator architecture for generative adversarial networks}.{\BBCQ}
\newblock
 \APACrefbtitle {{Proc. CVPR}} {{Proc. CVPR}}\ (\BPGS\ 4401--4410).
\PrintBackRefs{\CurrentBib}

\bibitem [\protect \citeauthoryear {%
Kingma%
\ \BBA {} Welling%
}{%
Kingma%
\ \BBA {} Welling%
}{%
{\protect \APACyear {2014}}%
}]{%
kingma2013auto}
\APACinsertmetastar {%
kingma2013auto}%
\begin{APACrefauthors}%
Kingma, D.P.%
\BCBT {}\ \BBA {} Welling, M.%
\end{APACrefauthors}%
\unskip\
\newblock
\APACrefYearMonthDay{2014}{}{}.
\newblock
{\BBOQ}\APACrefatitle {{Auto-encoding variational Bayes}} {{Auto-encoding variational Bayes}}.{\BBCQ}
\newblock
 \APACrefbtitle {{Proc. ICLR}} {{Proc. ICLR}}\ (\BPGS\ 1--14).
\PrintBackRefs{\CurrentBib}

\bibitem [\protect \citeauthoryear {%
Kohler%
, Krzyzak%
\BCBL {}\ \BBA {} Walter%
}{%
Kohler%
\ \protect \BOthers {.}}{%
{\protect \APACyear {2020}}%
}]{%
kohler2020rate}
\APACinsertmetastar {%
kohler2020rate}%
\begin{APACrefauthors}%
Kohler, M.%
, Krzyzak, A.%
\BCBL {} Walter, B.%
\end{APACrefauthors}%
\unskip\
\newblock
\APACrefYearMonthDay{2020}{}{}.
\newblock
{\BBOQ}\APACrefatitle {On the rate of convergence of image classifiers based on convolutional neural networks} {On the rate of convergence of image classifiers based on convolutional neural networks}.{\BBCQ}
\newblock
\APACjournalVolNumPages{ArXiv:2003.01526}{}{}{,}
\newblock

\newblock

\PrintBackRefs{\CurrentBib}

\bibitem [\protect \citeauthoryear {%
Kundu%
\ \BBA {} Dunson%
}{%
Kundu%
\ \BBA {} Dunson%
}{%
{\protect \APACyear {2014}}%
}]{%
kundu2014latent}
\APACinsertmetastar {%
kundu2014latent}%
\begin{APACrefauthors}%
Kundu, S.%
\BCBT {}\ \BBA {} Dunson, D.B.%
\end{APACrefauthors}%
\unskip\
\newblock
\APACrefYearMonthDay{2014}{}{}.
\newblock
{\BBOQ}\APACrefatitle {Latent factor models for density estimation} {Latent factor models for density estimation}.{\BBCQ}
\newblock
\APACjournalVolNumPages{Biometrika}{101}{3}{641--654,}
\newblock

\newblock

\PrintBackRefs{\CurrentBib}

\bibitem [\protect \citeauthoryear {%
Li%
, Chang%
, Cheng%
, Yang%
\BCBL {}\ \BBA {} P{\'o}czos%
}{%
Li%
\ \protect \BOthers {.}}{%
{\protect \APACyear {2017}}%
}]{%
li2017mmd}
\APACinsertmetastar {%
li2017mmd}%
\begin{APACrefauthors}%
Li, C\BHBI L.%
, Chang, W\BHBI C.%
, Cheng, Y.%
, Yang, Y.%
\BCBL {} P{\'o}czos, B.%
\end{APACrefauthors}%
\unskip\
\newblock
\APACrefYearMonthDay{2017}{}{}.
\newblock
{\BBOQ}\APACrefatitle {{MMD GAN: Towards deeper understanding of moment matching network}} {{MMD GAN: Towards deeper understanding of moment matching network}}.{\BBCQ}
\newblock
 \APACrefbtitle {{Proc. NIPS}} {{Proc. NIPS}}\ (\BPGS\ 2203--2213).
\PrintBackRefs{\CurrentBib}

\bibitem [\protect \citeauthoryear {%
Liang%
}{%
Liang%
}{%
{\protect \APACyear {2021}}%
}]{%
liang2021well}
\APACinsertmetastar {%
liang2021well}%
\begin{APACrefauthors}%
Liang, T.%
\end{APACrefauthors}%
\unskip\
\newblock
\APACrefYearMonthDay{2021}{}{}.
\newblock
{\BBOQ}\APACrefatitle {How well generative adversarial networks learn distributions} {How well generative adversarial networks learn distributions}.{\BBCQ}
\newblock
\APACjournalVolNumPages{J. Mach. Learn. Res.}{22}{228}{1--41,}
\newblock

\newblock

\PrintBackRefs{\CurrentBib}

\bibitem [\protect \citeauthoryear {%
Liu%
, Bousquet%
\BCBL {}\ \BBA {} Chaudhuri%
}{%
Liu%
\ \protect \BOthers {.}}{%
{\protect \APACyear {2017}}%
}]{%
liu2017approximation}
\APACinsertmetastar {%
liu2017approximation}%
\begin{APACrefauthors}%
Liu, S.%
, Bousquet, O.%
\BCBL {} Chaudhuri, K.%
\end{APACrefauthors}%
\unskip\
\newblock
\APACrefYearMonthDay{2017}{}{}.
\newblock
{\BBOQ}\APACrefatitle {Approximation and convergence properties of generative adversarial learning} {Approximation and convergence properties of generative adversarial learning}.{\BBCQ}
\newblock
 \APACrefbtitle {{Proc. NIPS}} {{Proc. NIPS}}\ (\BPGS\ 5545--5553).
\PrintBackRefs{\CurrentBib}

\bibitem [\protect \citeauthoryear {%
Manole%
, Balakrishnan%
, Niles-Weed%
\BCBL {}\ \BBA {} Wasserman%
}{%
Manole%
\ \protect \BOthers {.}}{%
{\protect \APACyear {2021}}%
}]{%
manole2021plugin}
\APACinsertmetastar {%
manole2021plugin}%
\begin{APACrefauthors}%
Manole, T.%
, Balakrishnan, S.%
, Niles-Weed, J.%
\BCBL {} Wasserman, L.%
\end{APACrefauthors}%
\unskip\
\newblock
\APACrefYearMonthDay{2021}{}{}.
\newblock
{\BBOQ}\APACrefatitle {Plugin estimation of smooth optimal transport maps} {Plugin estimation of smooth optimal transport maps}.{\BBCQ}
\newblock
\APACjournalVolNumPages{ArXiv:2107.12364}{}{}{,}
\newblock

\newblock

\PrintBackRefs{\CurrentBib}

\bibitem [\protect \citeauthoryear {%
Meister%
}{%
Meister%
}{%
{\protect \APACyear {2009}}%
}]{%
alexander2009deconvolution}
\APACinsertmetastar {%
alexander2009deconvolution}%
\begin{APACrefauthors}%
Meister, A.%
\end{APACrefauthors}%
\unskip\
\newblock
\APACrefYear{2009}.
\newblock
\APACrefbtitle {{Deconvolution Problems in Nonparametric Statistics}} {{Deconvolution Problems in Nonparametric Statistics}}.
\newblock
\APACaddressPublisher{}{Springer, New York}.
\PrintBackRefs{\CurrentBib}

\bibitem [\protect \citeauthoryear {%
Mescheder%
, Geiger%
\BCBL {}\ \BBA {} Nowozin%
}{%
Mescheder%
\ \protect \BOthers {.}}{%
{\protect \APACyear {2018}}%
}]{%
mescheder2018training}
\APACinsertmetastar {%
mescheder2018training}%
\begin{APACrefauthors}%
Mescheder, L.%
, Geiger, A.%
\BCBL {} Nowozin, S.%
\end{APACrefauthors}%
\unskip\
\newblock
\APACrefYearMonthDay{2018}{}{}.
\newblock
{\BBOQ}\APACrefatitle {{Which training methods for GANs do actually converge?}} {{Which training methods for GANs do actually converge?}}{\BBCQ}
\newblock
 \APACrefbtitle {{Proc. ICML}} {{Proc. ICML}}\ (\BPGS\ 3481--3490).
\PrintBackRefs{\CurrentBib}

\bibitem [\protect \citeauthoryear {%
Mroueh%
, Li%
, Sercu%
, Raj%
\BCBL {}\ \BBA {} Cheng%
}{%
Mroueh%
\ \protect \BOthers {.}}{%
{\protect \APACyear {2017}}%
}]{%
mroueh2017sobolev}
\APACinsertmetastar {%
mroueh2017sobolev}%
\begin{APACrefauthors}%
Mroueh, Y.%
, Li, C\BHBI L.%
, Sercu, T.%
, Raj, A.%
\BCBL {} Cheng, Y.%
\end{APACrefauthors}%
\unskip\
\newblock
\APACrefYearMonthDay{2017}{}{}.
\newblock
{\BBOQ}\APACrefatitle {{Sobolev GAN}} {{Sobolev GAN}}.{\BBCQ}
\newblock
\APACjournalVolNumPages{ArXiv:1711.04894}{}{}{,}
\newblock

\newblock

\PrintBackRefs{\CurrentBib}

\bibitem [\protect \citeauthoryear {%
M{\"u}ller%
}{%
M{\"u}ller%
}{%
{\protect \APACyear {1997}}%
}]{%
muller1997integral}
\APACinsertmetastar {%
muller1997integral}%
\begin{APACrefauthors}%
M{\"u}ller, A.%
\end{APACrefauthors}%
\unskip\
\newblock
\APACrefYearMonthDay{1997}{}{}.
\newblock
{\BBOQ}\APACrefatitle {Integral probability metrics and their generating classes of functions} {Integral probability metrics and their generating classes of functions}.{\BBCQ}
\newblock
\APACjournalVolNumPages{Adv. in Appl. Probab.}{29}{2}{429--443,}
\newblock

\newblock

\PrintBackRefs{\CurrentBib}

\bibitem [\protect \citeauthoryear {%
Nguyen%
}{%
Nguyen%
}{%
{\protect \APACyear {2013}}%
}]{%
nguyen2013convergence}
\APACinsertmetastar {%
nguyen2013convergence}%
\begin{APACrefauthors}%
Nguyen, X.%
\end{APACrefauthors}%
\unskip\
\newblock
\APACrefYearMonthDay{2013}{}{}.
\newblock
{\BBOQ}\APACrefatitle {Convergence of latent mixing measures in finite and infinite mixture models} {Convergence of latent mixing measures in finite and infinite mixture models}.{\BBCQ}
\newblock
\APACjournalVolNumPages{Ann. Statist.}{41}{1}{370--400,}
\newblock

\newblock

\PrintBackRefs{\CurrentBib}

\bibitem [\protect \citeauthoryear {%
Niles-Weed%
\ \BBA {} Berthet%
}{%
Niles-Weed%
\ \BBA {} Berthet%
}{%
{\protect \APACyear {2022}}%
}]{%
niles2022minimax}
\APACinsertmetastar {%
niles2022minimax}%
\begin{APACrefauthors}%
Niles-Weed, J.%
\BCBT {}\ \BBA {} Berthet, Q.%
\end{APACrefauthors}%
\unskip\
\newblock
\APACrefYearMonthDay{2022}{}{}.
\newblock
{\BBOQ}\APACrefatitle {{Minimax estimation of smooth densities in Wasserstein distance}} {{Minimax estimation of smooth densities in Wasserstein distance}}.{\BBCQ}
\newblock
\APACjournalVolNumPages{Ann. Statist.}{50}{3}{1519--1540,}
\newblock

\newblock

\PrintBackRefs{\CurrentBib}

\bibitem [\protect \citeauthoryear {%
Ohn%
\ \BBA {} Kim%
}{%
Ohn%
\ \BBA {} Kim%
}{%
{\protect \APACyear {2019}}%
}]{%
ohn2019smooth}
\APACinsertmetastar {%
ohn2019smooth}%
\begin{APACrefauthors}%
Ohn, I.%
\BCBT {}\ \BBA {} Kim, Y.%
\end{APACrefauthors}%
\unskip\
\newblock
\APACrefYearMonthDay{2019}{}{}.
\newblock
{\BBOQ}\APACrefatitle {Smooth function approximation by deep neural networks with general activation functions} {Smooth function approximation by deep neural networks with general activation functions}.{\BBCQ}
\newblock
\APACjournalVolNumPages{Entropy}{21}{7}{627,}
\newblock

\newblock

\PrintBackRefs{\CurrentBib}

\bibitem [\protect \citeauthoryear {%
Pati%
, Bhattacharya%
\BCBL {}\ \BBA {} Dunson%
}{%
Pati%
\ \protect \BOthers {.}}{%
{\protect \APACyear {2011}}%
}]{%
pati2011posterior}
\APACinsertmetastar {%
pati2011posterior}%
\begin{APACrefauthors}%
Pati, D.%
, Bhattacharya, A.%
\BCBL {} Dunson, D.B.%
\end{APACrefauthors}%
\unskip\
\newblock
\APACrefYearMonthDay{2011}{}{}.
\newblock
{\BBOQ}\APACrefatitle {Posterior convergence rates in non-linear latent variable models} {Posterior convergence rates in non-linear latent variable models}.{\BBCQ}
\newblock
\APACjournalVolNumPages{ArXiv:1109.5000}{}{}{,}
\newblock

\newblock

\PrintBackRefs{\CurrentBib}

\bibitem [\protect \citeauthoryear {%
Pooladian%
\ \BBA {} Niles-Weed%
}{%
Pooladian%
\ \BBA {} Niles-Weed%
}{%
{\protect \APACyear {2021}}%
}]{%
pooladian2021entropic}
\APACinsertmetastar {%
pooladian2021entropic}%
\begin{APACrefauthors}%
Pooladian, A\BHBI A.%
\BCBT {}\ \BBA {} Niles-Weed, J.%
\end{APACrefauthors}%
\unskip\
\newblock
\APACrefYearMonthDay{2021}{}{}.
\newblock
{\BBOQ}\APACrefatitle {Entropic estimation of optimal transport maps} {Entropic estimation of optimal transport maps}.{\BBCQ}
\newblock
\APACjournalVolNumPages{ArXiv:2109.12004}{}{}{,}
\newblock

\newblock

\PrintBackRefs{\CurrentBib}

\bibitem [\protect \citeauthoryear {%
Puchkin%
\ \BBA {} Spokoiny%
}{%
Puchkin%
\ \BBA {} Spokoiny%
}{%
{\protect \APACyear {2022}}%
}]{%
puchkin2022structure}
\APACinsertmetastar {%
puchkin2022structure}%
\begin{APACrefauthors}%
Puchkin, N.%
\BCBT {}\ \BBA {} Spokoiny, V.G.%
\end{APACrefauthors}%
\unskip\
\newblock
\APACrefYearMonthDay{2022}{}{}.
\newblock
{\BBOQ}\APACrefatitle {Structure-adaptive manifold estimation.} {Structure-adaptive manifold estimation.}{\BBCQ}
\newblock
\APACjournalVolNumPages{J. Mach. Learn. Res.}{23}{}{1--62,}
\newblock

\newblock

\PrintBackRefs{\CurrentBib}

\bibitem [\protect \citeauthoryear {%
Radford%
, Metz%
\BCBL {}\ \BBA {} Chintala%
}{%
Radford%
\ \protect \BOthers {.}}{%
{\protect \APACyear {2016}}%
}]{%
radford2016unsupervised}
\APACinsertmetastar {%
radford2016unsupervised}%
\begin{APACrefauthors}%
Radford, A.%
, Metz, L.%
\BCBL {} Chintala, S.%
\end{APACrefauthors}%
\unskip\
\newblock
\APACrefYearMonthDay{2016}{}{}.
\newblock
{\BBOQ}\APACrefatitle {Unsupervised representation learning with deep convolutional generative adversarial networks} {Unsupervised representation learning with deep convolutional generative adversarial networks}.{\BBCQ}
\newblock
 \APACrefbtitle {{Proc. ICLR}} {{Proc. ICLR}}\ (\BPGS\ 1--16).
\PrintBackRefs{\CurrentBib}

\bibitem [\protect \citeauthoryear {%
Rezende%
, Mohamed%
\BCBL {}\ \BBA {} Wierstra%
}{%
Rezende%
\ \protect \BOthers {.}}{%
{\protect \APACyear {2014}}%
}]{%
rezende2014stochastic}
\APACinsertmetastar {%
rezende2014stochastic}%
\begin{APACrefauthors}%
Rezende, D.J.%
, Mohamed, S.%
\BCBL {} Wierstra, D.%
\end{APACrefauthors}%
\unskip\
\newblock
\APACrefYearMonthDay{2014}{}{}.
\newblock
{\BBOQ}\APACrefatitle {Stochastic backpropagation and approximate inference in deep generative models} {Stochastic backpropagation and approximate inference in deep generative models}.{\BBCQ}
\newblock
 \APACrefbtitle {{Proc. ICML}} {{Proc. ICML}}\ (\BPGS\ 1278--1286).
\PrintBackRefs{\CurrentBib}

\bibitem [\protect \citeauthoryear {%
Schmidt-Hieber%
}{%
Schmidt-Hieber%
}{%
{\protect \APACyear {2020}}%
}]{%
schmidt2020nonparametric}
\APACinsertmetastar {%
schmidt2020nonparametric}%
\begin{APACrefauthors}%
Schmidt-Hieber, J.%
\end{APACrefauthors}%
\unskip\
\newblock
\APACrefYearMonthDay{2020}{}{}.
\newblock
{\BBOQ}\APACrefatitle {{Nonparametric regression using deep neural networks with ReLU activation function}} {{Nonparametric regression using deep neural networks with ReLU activation function}}.{\BBCQ}
\newblock
\APACjournalVolNumPages{Ann. Statist.}{48}{4}{1875--1897,}
\newblock

\newblock

\PrintBackRefs{\CurrentBib}

\bibitem [\protect \citeauthoryear {%
Schreuder%
}{%
Schreuder%
}{%
{\protect \APACyear {2021}}%
}]{%
schreuder2021bounding}
\APACinsertmetastar {%
schreuder2021bounding}%
\begin{APACrefauthors}%
Schreuder, N.%
\end{APACrefauthors}%
\unskip\
\newblock
\APACrefYearMonthDay{2021}{}{}.
\newblock
{\BBOQ}\APACrefatitle {{Bounding the expectation of the supremum of empirical processes indexed by H\"{o}lder classes}} {{Bounding the expectation of the supremum of empirical processes indexed by H\"{o}lder classes}}.{\BBCQ}
\newblock
\APACjournalVolNumPages{Math. Methods Statist.}{29}{}{76--86,}
\newblock

\newblock

\PrintBackRefs{\CurrentBib}

\bibitem [\protect \citeauthoryear {%
Schreuder%
, Brunel%
\BCBL {}\ \BBA {} Dalalyan%
}{%
Schreuder%
\ \protect \BOthers {.}}{%
{\protect \APACyear {2021}}%
}]{%
schreuder2021statistical}
\APACinsertmetastar {%
schreuder2021statistical}%
\begin{APACrefauthors}%
Schreuder, N.%
, Brunel, V\BHBI E.%
\BCBL {} Dalalyan, A.%
\end{APACrefauthors}%
\unskip\
\newblock
\APACrefYearMonthDay{2021}{}{}.
\newblock
{\BBOQ}\APACrefatitle {Statistical guarantees for generative models without domination} {Statistical guarantees for generative models without domination}.{\BBCQ}
\newblock
 \APACrefbtitle {{Proc. Algorithmic Learning Theory}} {{Proc. Algorithmic Learning Theory}}\ (\BPGS\ 1051--1071).
\PrintBackRefs{\CurrentBib}

\bibitem [\protect \citeauthoryear {%
Singh%
\ \BBA {} P{\'o}czos%
}{%
Singh%
\ \BBA {} P{\'o}czos%
}{%
{\protect \APACyear {2018}}%
}]{%
singh2018minimax}
\APACinsertmetastar {%
singh2018minimax}%
\begin{APACrefauthors}%
Singh, S.%
\BCBT {}\ \BBA {} P{\'o}czos, B.%
\end{APACrefauthors}%
\unskip\
\newblock
\APACrefYearMonthDay{2018}{}{}.
\newblock
{\BBOQ}\APACrefatitle {{Minimax distribution estimation in Wasserstein distance}} {{Minimax distribution estimation in Wasserstein distance}}.{\BBCQ}
\newblock
\APACjournalVolNumPages{ArXiv:1802.08855}{}{}{,}
\newblock

\newblock

\PrintBackRefs{\CurrentBib}

\bibitem [\protect \citeauthoryear {%
Singh%
\ \protect \BOthers {.}}{%
Singh%
\ \protect \BOthers {.}}{%
{\protect \APACyear {2018}}%
}]{%
singh2018nonparametric}
\APACinsertmetastar {%
singh2018nonparametric}%
\begin{APACrefauthors}%
Singh, S.%
, Uppal, A.%
, Li, B.%
, Li, C\BHBI L.%
, Zaheer, M.%
\BCBL {} P{\'o}czos, B.%
\end{APACrefauthors}%
\unskip\
\newblock
\APACrefYearMonthDay{2018}{}{}.
\newblock
{\BBOQ}\APACrefatitle {Nonparametric density estimation with adversarial losses} {Nonparametric density estimation with adversarial losses}.{\BBCQ}
\newblock
 \APACrefbtitle {{Proc. NeurIPS}} {{Proc. NeurIPS}}\ (\BPGS\ 10246--10257).
\PrintBackRefs{\CurrentBib}

\bibitem [\protect \citeauthoryear {%
St{\'e}phanovitch%
, Aamari%
\BCBL {}\ \BBA {} Levrard%
}{%
St{\'e}phanovitch%
\ \protect \BOthers {.}}{%
{\protect \APACyear {2023}}%
}]{%
stephanovitch2023wasserstein}
\APACinsertmetastar {%
stephanovitch2023wasserstein}%
\begin{APACrefauthors}%
St{\'e}phanovitch, A.%
, Aamari, E.%
\BCBL {} Levrard, C.%
\end{APACrefauthors}%
\unskip\
\newblock
\APACrefYearMonthDay{2023}{}{}.
\newblock
{\BBOQ}\APACrefatitle {Wasserstein generative adversarial networks are minimax optimal distribution estimators} {Wasserstein generative adversarial networks are minimax optimal distribution estimators}.{\BBCQ}
\newblock
\APACjournalVolNumPages{ArXiv:2311.18613}{}{}{,}
\newblock

\newblock

\PrintBackRefs{\CurrentBib}

\bibitem [\protect \citeauthoryear {%
Tang%
\ \BBA {} Yang%
}{%
Tang%
\ \BBA {} Yang%
}{%
{\protect \APACyear {2023}}%
}]{%
tang2023minimax}
\APACinsertmetastar {%
tang2023minimax}%
\begin{APACrefauthors}%
Tang, R.%
\BCBT {}\ \BBA {} Yang, Y.%
\end{APACrefauthors}%
\unskip\
\newblock
\APACrefYearMonthDay{2023}{}{}.
\newblock
{\BBOQ}\APACrefatitle {Minimax rate of distribution estimation on unknown submanifolds under adversarial losses} {Minimax rate of distribution estimation on unknown submanifolds under adversarial losses}.{\BBCQ}
\newblock
\APACjournalVolNumPages{Ann. Statist.}{51}{3}{1282--1308,}
\newblock

\newblock

\PrintBackRefs{\CurrentBib}

\bibitem [\protect \citeauthoryear {%
Telgarsky%
}{%
Telgarsky%
}{%
{\protect \APACyear {2016}}%
}]{%
telgarsky2016benefits}
\APACinsertmetastar {%
telgarsky2016benefits}%
\begin{APACrefauthors}%
Telgarsky, M.%
\end{APACrefauthors}%
\unskip\
\newblock
\APACrefYearMonthDay{2016}{}{}.
\newblock
{\BBOQ}\APACrefatitle {Benefits of depth in neural networks} {Benefits of depth in neural networks}.{\BBCQ}
\newblock
 \APACrefbtitle {{Proc. COLT}} {{Proc. COLT}}\ (\BPGS\ 1517--1539).
\PrintBackRefs{\CurrentBib}

\bibitem [\protect \citeauthoryear {%
Tsybakov%
}{%
Tsybakov%
}{%
{\protect \APACyear {2008}}%
}]{%
tsybakov2008introduction}
\APACinsertmetastar {%
tsybakov2008introduction}%
\begin{APACrefauthors}%
Tsybakov, A.B.%
\end{APACrefauthors}%
\unskip\
\newblock
\APACrefYear{2008}.
\newblock
\APACrefbtitle {{Introduction to Nonparametric Estimation}} {{Introduction to Nonparametric Estimation}}.
\newblock
\APACaddressPublisher{}{Springer, New York}.
\PrintBackRefs{\CurrentBib}

\bibitem [\protect \citeauthoryear {%
Uppal%
, Singh%
\BCBL {}\ \BBA {} P{\'o}czos%
}{%
Uppal%
\ \protect \BOthers {.}}{%
{\protect \APACyear {2019}}%
}]{%
uppal2019nonparametric}
\APACinsertmetastar {%
uppal2019nonparametric}%
\begin{APACrefauthors}%
Uppal, A.%
, Singh, S.%
\BCBL {} P{\'o}czos, B.%
\end{APACrefauthors}%
\unskip\
\newblock
\APACrefYearMonthDay{2019}{}{}.
\newblock
{\BBOQ}\APACrefatitle {{Nonparametric density estimation and convergence of GANs under Besov IPM losses}} {{Nonparametric density estimation and convergence of GANs under Besov IPM losses}}.{\BBCQ}
\newblock
 \APACrefbtitle {{Proc. NeurIPS}} {{Proc. NeurIPS}}\ (\BPGS\ 9089--9100).
\PrintBackRefs{\CurrentBib}

\bibitem [\protect \citeauthoryear {%
Urbas%
}{%
Urbas%
}{%
{\protect \APACyear {1988}}%
}]{%
urbas1988regularity}
\APACinsertmetastar {%
urbas1988regularity}%
\begin{APACrefauthors}%
Urbas, J.I.%
\end{APACrefauthors}%
\unskip\
\newblock
\APACrefYearMonthDay{1988}{}{}.
\newblock
{\BBOQ}\APACrefatitle {{Regularity of generalized solutions of Monge--Amp\`{e}re equations}} {{Regularity of generalized solutions of Monge--Amp\`{e}re equations}}.{\BBCQ}
\newblock
\APACjournalVolNumPages{Math. Z.}{197}{3}{365--393,}
\newblock

\newblock

\PrintBackRefs{\CurrentBib}

\bibitem [\protect \citeauthoryear {%
van~der Vaart%
\ \BBA {} Wellner%
}{%
van~der Vaart%
\ \BBA {} Wellner%
}{%
{\protect \APACyear {1996}}%
}]{%
van1996weak}
\APACinsertmetastar {%
van1996weak}%
\begin{APACrefauthors}%
van~der Vaart, A.W.%
\BCBT {}\ \BBA {} Wellner, J.A.%
\end{APACrefauthors}%
\unskip\
\newblock
\APACrefYear{1996}.
\newblock
\APACrefbtitle {{Weak Convergence and Empirical Processes}} {{Weak Convergence and Empirical Processes}}.
\newblock
\APACaddressPublisher{}{Springer}.
\PrintBackRefs{\CurrentBib}

\bibitem [\protect \citeauthoryear {%
Villani%
}{%
Villani%
}{%
{\protect \APACyear {2003}}%
}]{%
villani2003topics}
\APACinsertmetastar {%
villani2003topics}%
\begin{APACrefauthors}%
Villani, C.%
\end{APACrefauthors}%
\unskip\
\newblock
\APACrefYear{2003}.
\newblock
\APACrefbtitle {{Topics in Optimal Transportation}} {{Topics in Optimal Transportation}}.
\newblock
\APACaddressPublisher{}{American Mathematical Society}.
\PrintBackRefs{\CurrentBib}

\bibitem [\protect \citeauthoryear {%
Villani%
}{%
Villani%
}{%
{\protect \APACyear {2008}}%
}]{%
villani2008optimal}
\APACinsertmetastar {%
villani2008optimal}%
\begin{APACrefauthors}%
Villani, C.%
\end{APACrefauthors}%
\unskip\
\newblock
\APACrefYear{2008}.
\newblock
\APACrefbtitle {{Optimal Transport: Old and New}} {{Optimal Transport: Old and New}}.
\newblock
\APACaddressPublisher{}{Springer}.
\PrintBackRefs{\CurrentBib}

\bibitem [\protect \citeauthoryear {%
Wainwright%
}{%
Wainwright%
}{%
{\protect \APACyear {2019}}%
}]{%
wainwright2019high}
\APACinsertmetastar {%
wainwright2019high}%
\begin{APACrefauthors}%
Wainwright, M.J.%
\end{APACrefauthors}%
\unskip\
\newblock
\APACrefYear{2019}.
\newblock
\APACrefbtitle {{High-Dimensional Statistics: A Non-Asymptotic Viewpoint}} {{High-Dimensional Statistics: A Non-Asymptotic Viewpoint}}.
\newblock
\APACaddressPublisher{}{Cambridge University Press}.
\PrintBackRefs{\CurrentBib}

\bibitem [\protect \citeauthoryear {%
Weed%
\ \BBA {} Bach%
}{%
Weed%
\ \BBA {} Bach%
}{%
{\protect \APACyear {2019}}%
}]{%
weed2019sharp}
\APACinsertmetastar {%
weed2019sharp}%
\begin{APACrefauthors}%
Weed, J.%
\BCBT {}\ \BBA {} Bach, F.%
\end{APACrefauthors}%
\unskip\
\newblock
\APACrefYearMonthDay{2019}{}{}.
\newblock
{\BBOQ}\APACrefatitle {{Sharp asymptotic and finite-sample rates of convergence of empirical measures in Wasserstein distance}} {{Sharp asymptotic and finite-sample rates of convergence of empirical measures in Wasserstein distance}}.{\BBCQ}
\newblock
\APACjournalVolNumPages{Bernoulli}{25}{4A}{2620--2648,}
\newblock

\newblock

\PrintBackRefs{\CurrentBib}

\bibitem [\protect \citeauthoryear {%
Wei%
\ \BBA {} Nguyen%
}{%
Wei%
\ \BBA {} Nguyen%
}{%
{\protect \APACyear {2022}}%
}]{%
wei2022convergence}
\APACinsertmetastar {%
wei2022convergence}%
\begin{APACrefauthors}%
Wei, Y.%
\BCBT {}\ \BBA {} Nguyen, X.%
\end{APACrefauthors}%
\unskip\
\newblock
\APACrefYearMonthDay{2022}{}{}.
\newblock
{\BBOQ}\APACrefatitle {{Convergence of de Finetti's mixing measure in latent structure models for observed exchangeable sequences}} {{Convergence of de Finetti's mixing measure in latent structure models for observed exchangeable sequences}}.{\BBCQ}
\newblock
\APACjournalVolNumPages{Ann. Statist.}{50}{4}{1859 -- 1889,}
\newblock

\newblock

\PrintBackRefs{\CurrentBib}

\bibitem [\protect \citeauthoryear {%
Wong%
\ \BBA {} Shen%
}{%
Wong%
\ \BBA {} Shen%
}{%
{\protect \APACyear {1995}}%
}]{%
wong1995probability}
\APACinsertmetastar {%
wong1995probability}%
\begin{APACrefauthors}%
Wong, W.H.%
\BCBT {}\ \BBA {} Shen, X.%
\end{APACrefauthors}%
\unskip\
\newblock
\APACrefYearMonthDay{1995}{}{}.
\newblock
{\BBOQ}\APACrefatitle {{Probability inequalities for likelihood ratios and convergence rates of sieve MLEs}} {{Probability inequalities for likelihood ratios and convergence rates of sieve MLEs}}.{\BBCQ}
\newblock
\APACjournalVolNumPages{Ann. Statist.}{23}{2}{339--362,}
\newblock

\newblock

\PrintBackRefs{\CurrentBib}

\bibitem [\protect \citeauthoryear {%
Yalcin%
\ \BBA {} Amemiya%
}{%
Yalcin%
\ \BBA {} Amemiya%
}{%
{\protect \APACyear {2001}}%
}]{%
yalcin2001nonlinear}
\APACinsertmetastar {%
yalcin2001nonlinear}%
\begin{APACrefauthors}%
Yalcin, I.%
\BCBT {}\ \BBA {} Amemiya, Y.%
\end{APACrefauthors}%
\unskip\
\newblock
\APACrefYearMonthDay{2001}{}{}.
\newblock
{\BBOQ}\APACrefatitle {Nonlinear factor analysis as a statistical method} {Nonlinear factor analysis as a statistical method}.{\BBCQ}
\newblock
\APACjournalVolNumPages{Statist. Sci.}{16}{3}{275--294,}
\newblock

\newblock

\PrintBackRefs{\CurrentBib}

\bibitem [\protect \citeauthoryear {%
Yarotsky%
}{%
Yarotsky%
}{%
{\protect \APACyear {2017}}%
}]{%
yarotsky2017error}
\APACinsertmetastar {%
yarotsky2017error}%
\begin{APACrefauthors}%
Yarotsky, D.%
\end{APACrefauthors}%
\unskip\
\newblock
\APACrefYearMonthDay{2017}{}{}.
\newblock
{\BBOQ}\APACrefatitle {{Error bounds for approximations with deep ReLU networks}} {{Error bounds for approximations with deep ReLU networks}}.{\BBCQ}
\newblock
\APACjournalVolNumPages{Neural Networks}{94}{}{103--114,}
\newblock

\newblock

\PrintBackRefs{\CurrentBib}

\bibitem [\protect \citeauthoryear {%
Yarotsky%
}{%
Yarotsky%
}{%
{\protect \APACyear {2021}}%
}]{%
yarotsky2021universal}
\APACinsertmetastar {%
yarotsky2021universal}%
\begin{APACrefauthors}%
Yarotsky, D.%
\end{APACrefauthors}%
\unskip\
\newblock
\APACrefYearMonthDay{2021}{}{}.
\newblock
{\BBOQ}\APACrefatitle {Universal approximations of invariant maps by neural networks} {Universal approximations of invariant maps by neural networks}.{\BBCQ}
\newblock
\APACjournalVolNumPages{Constr. Approx.}{}{}{1--68,}
\newblock

\newblock

\PrintBackRefs{\CurrentBib}

\bibitem [\protect \citeauthoryear {%
Zhang%
, Liu%
, Zhou%
, Xu%
\BCBL {}\ \BBA {} He%
}{%
Zhang%
\ \protect \BOthers {.}}{%
{\protect \APACyear {2018}}%
}]{%
zhang2018discrimination}
\APACinsertmetastar {%
zhang2018discrimination}%
\begin{APACrefauthors}%
Zhang, P.%
, Liu, Q.%
, Zhou, D.%
, Xu, T.%
\BCBL {} He, X.%
\end{APACrefauthors}%
\unskip\
\newblock
\APACrefYearMonthDay{2018}{}{}.
\newblock
{\BBOQ}\APACrefatitle {{On the discrimination-generalization tradeoff in GANs}} {{On the discrimination-generalization tradeoff in GANs}}.{\BBCQ}
\newblock
 \APACrefbtitle {{Proc. ICLR}} {{Proc. ICLR}}\ (\BPGS\ 1--26).
\PrintBackRefs{\CurrentBib}

\end{thebibliography}
\addtocontents{toc}{\protect\setcounter{tocdepth}{2}}

\pagebreak

\tableofcontents

\begin{appendix}

\section{Proof of Theorem \ref{thm:rate-general}}

Choose $\bg_* \in \cG$ such that
\be\label{eq:general-tech2}
	d_{\rm eval}(Q_*, Q_0) \leq \inf_{\bg\in\cG} d_{\rm eval} (Q_\bg, Q_0) + \epsilon_1 \stackrel{\rm (i)}{\leq} 2\epsilon_1,
\ee
where $Q_* = Q_{\bg_*}$.
Then,
\bean
	&& d_{\rm eval}(\hat Q, Q_0) 
	\leq d_{\rm eval}( \hat Q, Q_*) + d_{\rm eval}(Q_*, Q_0)
	\\
	&& \stackrel{\eqref{eq:general-tech2}}{\leq} d_{\rm eval} ( \hat Q, Q_*) + 2\epsilon_1
	\stackrel{\rm (iv)}{\leq} d_\cF( \hat Q, Q_*) + 2\epsilon_1 + \epsilon_4
	\\
	&& \leq d_\cF(\hat Q, \bbP_n) + d_\cF(\bbP_n, Q_*) + 2\epsilon_1 + \epsilon_4
	\\
	&& \stackrel{\rm (ii)}{\leq} \inf_{\bg \in \cG} d_\cF (Q_\bg, \bbP_n) + d_\cF(\bbP_n, Q_*) +  2\epsilon_1 + \epsilon_2 + \epsilon_4
	\\
	&& \leq \inf_{\bg \in \cG} d_\cF (Q_\bg, \bbP_n) + d_\cF(\bbP_n, P_0) + d_\cF(P_0, Q_0) + d_\cF(Q_0, Q_*) + 2\epsilon_1 + \epsilon_2 + \epsilon_4
	\\
	&& \leq \inf_{\bg \in \cG} d_\cF (Q_\bg, Q_0) + d_\cF(\bbP_n, Q_0) + d_\cF(\bbP_n, P_0) + d_\cF(P_0, Q_0) + d_\cF(Q_0, Q_*) 
	\\
	&& \qquad + 2\epsilon_1 + \epsilon_2 + \epsilon_4
	\\
	&& \stackrel{\rm (iv)}{\leq} \inf_{\bg \in \cG} d_{\rm eval} (Q_\bg, Q_0) + d_\cF(\bbP_n, Q_0) + d_\cF(\bbP_n, P_0) + d_\cF(P_0, Q_0) + d_\cF(Q_0, Q_*) 
	\\
	&& \qquad + 2\epsilon_1 + \epsilon_2 + 2\epsilon_4
	\\
	&& \stackrel{\rm (i)}{\leq}  d_\cF(\bbP_n, Q_0) + d_\cF(\bbP_n, P_0) + d_\cF(P_0, Q_0) + d_\cF(Q_0, Q_*) + 3\epsilon_1 + \epsilon_2 + 2\epsilon_4
	\\
	&& \leq 2 d_\cF(\bbP_n, P_0) + 2 d_\cF(P_0, Q_0) + d_\cF(Q_0, Q_*) + 3\epsilon_1 + \epsilon_2 + 2\epsilon_4
	\\
	&& \stackrel{\rm (iv)}{\leq} 2 d_\cF(\bbP_n, P_0) + 2 d_\cF(P_0, Q_0) + d_{\rm eval}(Q_0, Q_*) + 3\epsilon_1 + \epsilon_2 + 3\epsilon_4
	\\
	&& \stackrel{\eqref{eq:general-tech2}}{\leq} 2 d_\cF(\bbP_n, P_0) + 2 d_\cF(P_0, Q_0) + 5\epsilon_1 + \epsilon_2 + 3\epsilon_4.
\eean
By taking the expectation, we complete the proof.
\qed

\section{Proof of Theorem \ref{thm:rate-composition}}

We will construct a generator class $\cG$ and a discriminator $\cF$ satisfying condition \eqref{eq:general-condition} of Theorem \ref{thm:rate-general} with $d_{\rm eval} = W_1$.
By the construction of the estimator $\hat Q$, condition \eqref{eq:general-condition}-(ii) is automatically satisfied with $\epsilon_2 = \epsilon_{\rm opt}$ for any $\cG$ and $\cF$.

Let $\delta > 0$ be given.
Lemma 3.5 from \cite{chae2023likelihood} implies that there exists $\bg_* \in \cD(L, \bp, s, K\vee 1)$, with
\bean
	L \leq c_1 \log \delta^{-1}, 
	\quad |\bp|_\infty \leq c_1 \delta^{-t_*/\beta_*},
	\quad s \leq c_1 \delta^{-t_*/\beta_*} \log \delta^{-1}
\eean
for some constant $c_1 = c_1(q, \bd, \bt, \bbeta, K)$, such that $\| \bg_* - \bg_0 \|_\infty < \delta$.
Let $Q_* = Q_{\bg_*}$ and $\cG = \cD(L, \bp, s, K\vee 1)$.
Then, by the Kantorovich--Rubinstein duality (see Theorem 1.14 in \cite{villani2003topics}),
\bean
	&& W_1(Q_*, Q_0) 
	= \sup_{f \in \cF_{\rm Lip}} |Q_* f - Q_0 f|
	\\
	&& \leq \sup_{f \in \cF_{\rm Lip}} \int \Big| f\big(\bg_*(\bz)\big) - f\big(\bg_0(\bz)\big) \Big| dP_Z(\bz)
	\\
	&& \leq \int |\bg_*(\bz) - \bg_0(\bz)|_2 dP_Z(\bz)
	\leq \sqrt{D} \|\bg_* - \bg_0 \|_\infty 
	\leq \sqrt{D} \delta.
\eean
Hence, condition \eqref{eq:general-condition}-(i) holds with $\epsilon_1 = \sqrt{D}\delta$.

Let $\epsilon > 0$ be given.
For two Borel probability measures $Q_1$ and $Q_2$ on $\bbR^D$, one can choose $f_{Q_1, Q_2} \in \cF_{\rm Lip}$ such that $f_{Q_1, Q_2}(\bzero_D) = 0$ and
\bean
	W_1(Q_1, Q_2) 
	= \sup_{f\in\cF_{\rm Lip}} | Q_1 f - Q_2 f |
	\leq | Q_1 f_{Q_1, Q_2} - Q_2 f_{Q_1, Q_2}| + \epsilon.
\eean
Then, by the Lipschitz continuity, 
\bean
	\sup_{|\bx|_\infty \leq K} | f_{Q_1, Q_2}(\bx) | 
	\leq \sup_{|\bx|_\infty \leq K} |\bx|_2
	= \sqrt{D} K.
\eean
Let $\{\bg_1, \ldots, \bg_N\}$ be an $\epsilon$-cover of $\cG \cup \{\bg_0\}$ with respect to $\|\cdot\|_{P_Z, 2}$ and 
\bean
	\cF = \Big\{ f_{jk}: 1 \leq j, k \leq N \Big\},
\eean
where
\bean
	\|\bg\|_{P_Z, p} &=& \left(\int |\bg(\bz) |_p^p dP_Z(\bz) \right)^{1/p}
\eean
and $f_{jk} = f_{Q_{\bg_j}, Q_{\bg_k}}$.
Since $\|\bg - \widetilde \bg\|_{P_Z, 2} \leq \sqrt{D} \| \bg - \widetilde \bg \|_\infty$ for every $\bg, \widetilde \bg \in \cG \cup \{\bg_0\}$ and
\bean
	&& \log N(\epsilon, \cG, \| \cdot \|_\infty) 
	\\
	&& \leq (s+1) \bigg\{\log 2 + \log \epsilon^{-1} + \log (L+1) + 2\sum_{l=0}^{L+1} \log (p_l + 1) \bigg\}
\eean
by Lemma 5 of \cite{schmidt2020nonparametric}, the number $N$ can be bounded as
\be\begin{split}\label{eq:logN-bound}
	\log N 
	&\leq \log \big(N ( \epsilon/\sqrt{D}, \cG, \| \cdot\|_\infty ) +1\big)
	\\
	& \leq c_2 s \Big( \log D + \log \epsilon^{-1} + L \log \delta^{-1} \Big)
	\\
	&\leq c_3 \delta^{-t_*/\beta_*} \log \delta^{-1} \Big\{ \log \epsilon^{-1} + (\log \delta^{-1})^2 \Big\},
\end{split}\ee
where $c_2 = c_2(t_*, \beta_*)$ and $c_3 = c_3(c_1, c_2, D)$.
Here, $N(\epsilon, \cG, \|\cdot\|_\infty)$ denotes the covering number of $\cG$ with respect to $\|\cdot\|_\infty$.

Next, we will prove that condition \eqref{eq:general-condition}-(iv) is satisfied with $\epsilon_4 = 5\epsilon$.
Note that $d_\cF \leq W_1$ by the construction.
For $\bg, \widetilde\bg \in \cG \cup \{\bg_0\}$, we can choose $\bg_j$ and $\bg_k$ such that $\|\bg - \bg_j\|_{P_Z,2} \leq \epsilon$ and $\|\widetilde\bg - \bg_k\|_{P_Z,2} \leq \epsilon$.
Then,
\be\begin{split}\label{eq:comp-tech1}
	W_1(Q_\bg, Q_{\widetilde\bg}) 
	&\leq W_1(Q_\bg, Q_{\bg_j}) + W_1(Q_{\bg_j}, Q_{\bg_k}) + W_1(Q_{\bg_k}, Q_{\widetilde\bg})
	\\
	&\leq W_1(Q_\bg, Q_{\bg_j}) + d_\cF(Q_{\bg_j}, Q_{\bg_k}) + W_1(Q_{\bg_k}, Q_{\widetilde\bg}) + \epsilon.
\end{split}\ee
Note that
\bean
	W_1(Q_\bg, Q_{\bg_j}) 
	&=& \sup_{f\in\cF_{\rm Lip}} \left| \int f\big(\bg(\bz)\big) dP_Z(\bz) - \int f\big(\bg_j(\bz)\big) dP_Z(\bz) \right|
	\\
	&\leq& \int |\bg(\bz) - \bg_j(\bz)|_2 dP_Z(\bz)
	\leq \|\bg - \bg_j\|_{P_Z, 2} \leq \epsilon.
\eean
Similarly, $W_1(Q_{\bg_k}, Q_{\widetilde\bg}) \leq \epsilon$, and therefore,
\bean
	d_\cF(Q_{\bg_j}, Q_{\bg_k}) 
	&\leq& d_\cF(Q_{\bg_j}, Q_\bg) + d_\cF(Q_\bg, Q_{\widetilde\bg}) + d_\cF(Q_{\widetilde\bg}, Q_{\bg_k})
	\\
	&\leq& d_\cF(Q_\bg, Q_{\widetilde\bg}) + 2\epsilon.
\eean
Hence, the right hand side of \eqref{eq:comp-tech1} is bounded by $d_\cF(Q_\bg, Q_{\widetilde\bg}) + 5\epsilon$.
That is, condition \eqref{eq:general-condition}-(iv) holds with $\epsilon_4=5\epsilon$.

Next, note that $\bbP_n$ is the empirical measure based on \iid\ samples from $P_0$.
Let $\bY$ and $\bepsilon$ be independent random vectors following $Q_0$ and $\cN(\bzero_D, \sigma_0^2 \Id_D)$, respectively.
For any $f \in \cF$, by the Lipschitz continuity,
\bean
	|f(\bY + \bepsilon) | \leq |\bY + \bepsilon|_2 \leq |\bY|_2 + |\bepsilon|_2.
\eean
Since $\bY$ is bounded almost surely and $\sigma_0 \leq 1$, $f(\bY + \bepsilon)$ is a sub-Gaussian random variable with the sub-Gaussian parameter $\sigma=\sigma(K, D)$.
By the Hoeffding's inequality,
\bean
	P_0 \Big( \big|\bbP_n f - P_0 f \big| > t \Big) \leq 2 \exp\left[- \frac{nt^2}{2\sigma^2} \right]
\eean
for every $f \in \cF$ and $t \geq 0$; see Proposition 2.5 from \cite{wainwright2019high} for Hoeffding's inequality for unbounded sub-Gaussian random variables.
Since $\cF$ is a finite set with the cardinality $N^2$,
\bean
	P_0 \bigg( \sup_{f\in\cF}\big|\bbP_n f - P_0 f \big| > t \bigg) \leq 2N^2 \exp\left[- \frac{nt^2}{2\sigma^2} \right].
\eean
If $t \geq 2\sigma \sqrt{\{\log (2N^2)\}/n}$, the right hand side is bounded by $e^{-nt^2/(4\sigma^2)}$.
Therefore,
\bean
	&& \E d_\cF(\bbP_n, P_0)
	= \int_0^\infty P_0 \big( d_\cF(\bbP_n, P_0) > t \big) dt
	\\
	&& \leq 2\sigma \sqrt{\frac{ \log (2N^2)}{n}} + \int_0^\infty \exp\left[- \frac{nt^2}{4\sigma^2} \right] dt
	\\
	&& \leq 2\sigma \sqrt{\frac{\log (2N^2)}{n}} + \sigma\sqrt{\frac{\pi}{n}}
\eean
and condition \eqref{eq:general-condition}-(iii) is also satisfied with $\epsilon_3$ equal to the right hand side of the last display.

Note that
\bean
	d_\cF(P_0, Q_0) \leq W_1(P_0, Q_0) \leq W_2(P_0, Q_0) \leq \sqrt{D} \sigma_0,
\eean
where the last inequality holds because $P_0$ is the convolution of $Q_0$ and $\cN(\bzero_D, \sigma_0^2 \Id_D)$.
By Theorem \ref{thm:rate-general}, we have
\bean
	\E W_1(\hat Q, Q_0) 
	&\leq& 2 \sqrt{D} \sigma_0 + 5 \sqrt{D} \delta + \epsilon_{\rm opt} + 4 \sigma \sqrt{\frac{\log (2N^2)}{n}} + 2\sigma\sqrt{\frac{\pi}{n}} + 10\epsilon
	\\
	&\leq& c_4\bigg\{ \epsilon_{\rm opt} + \sigma_0 + \delta + \sqrt{\frac{\log N}{n}} + \epsilon \bigg\} ,
\eean
where $c_4 = c_4(\sigma, D)$.
Combining with \eqref{eq:logN-bound}, we have 
\bean
	\E W_1(\hat Q, Q_0)
	\leq c_5 \bigg\{ \epsilon_{\rm opt} + \sigma_0 + \delta + \frac{ \sqrt{\log \delta^{-1}} ( \sqrt{\log \epsilon^{-1}} + \log \delta^{-1})}{\sqrt{n} \delta^{t_*/2\beta_*}} + \epsilon \bigg\},
\eean
where $c_5 = c_5(c_3, c_4)$.
The proof is complete if we take 
\bean
	\delta = n^{-\beta_* / (2\beta_* + t_*)} ( \log n )^{\frac{3\beta_*}{2 \beta_* + t_*}}
\eean
and $\epsilon = n^{-\log n}$.
\qed

\section{Proof of Theorem \ref{thm:rate-ipm}}

The proof of the first assertion is the same as that of Theorem \ref{thm:rate-composition}.
The only difference is that some constants in the proof depend on the Lipschitz constant $C_1$.

For the second assertion, we utilize Theorem \ref{thm:rate-general} with $\cF = \cF_0$.
Since $d_{\rm eval} = d_\cF$, we have $\epsilon_4 = 0$.
Also, for a large enough $\cG$, \ie\ large depth, width and number of nonzero parameters, $\epsilon_1$ can be set to be an arbitrarily small number.
Since $\cF$ consists of Lipschitz continuous function, $d_\cF(P_0, Q_0) \lesssim \sigma_0$.
It follows by Theorem \ref{thm:rate-general} that $\E d_{\cF_0}(\hat Q, Q_0) \lesssim \sigma_0 + \epsilon_{\rm opt} + \E d_{\cF_0}(\bbP_n, P_0)$.

\section{Proof of Theorem \ref{thm:lower-bound}}

The proof is divided into several cases. For cases with $q=0$, we write $\beta_0$ as $\beta$ for simplicity.

\subsubsection*{Case 1: $q=0$ and $t_0 = d = D$}

In this case, $\beta_* = \beta$, $t_*=d$ and $\cG_0 = \cH_K^\beta ([0,1]^d) \times \cdots \times \cH_K^\beta ([0,1]^d)$.
Our proof relies on Fano's method for which we refer to Chapter 15 from \cite{wainwright2019high}.

\begin{figure}
	\begin{center}
		\subfigure[$\phi(z)$]{\scalebox{0.15}{\includegraphics{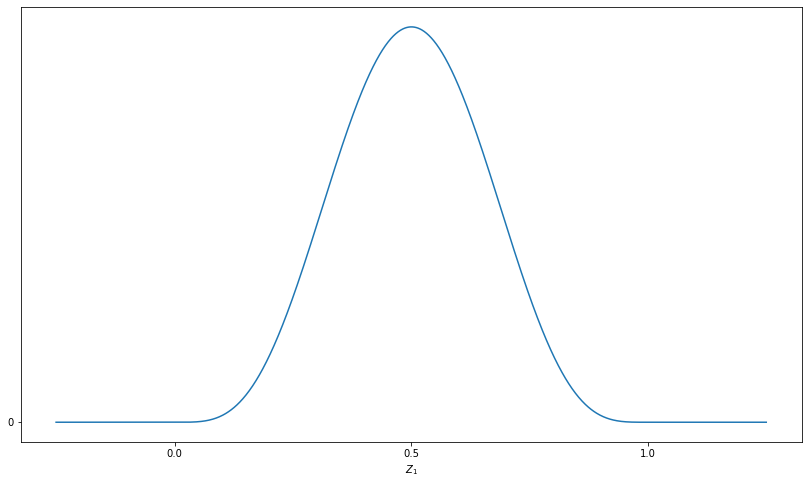}}}
		\subfigure[$\phi'(z)$]{\scalebox{0.15}{\includegraphics{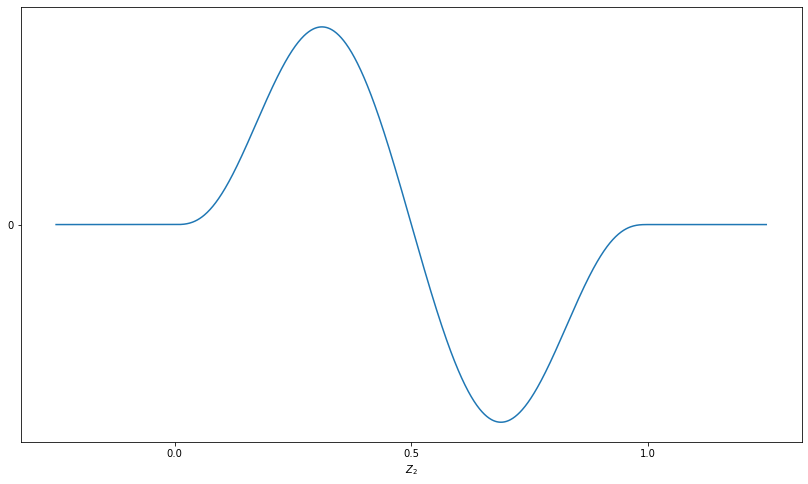}}}
		\subfigure[$\phi'(z_1) \phi(z_2)$]{\scalebox{0.042}{\includegraphics{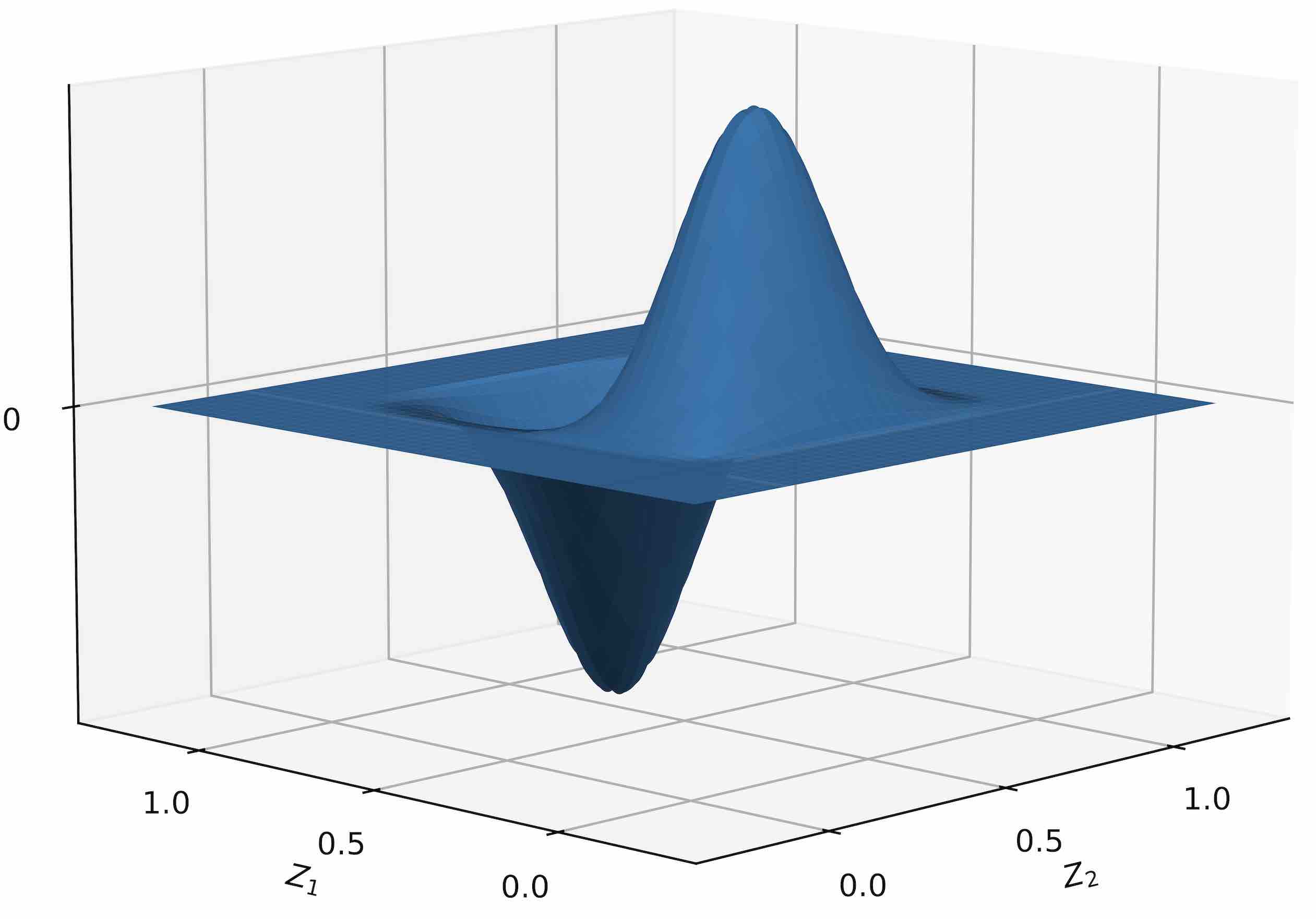}}}
		\caption{An illustration of $\phi$ and related functions.} \label{fig:phi}
	\end{center}
\end{figure}

Let $\phi: \bbR \to [0, \infty)$ be a fixed function satisfying that 
\ben
	\item[(i)] $\phi$ is $[\beta+1]$-times continuously differentiable on $\bbR$,
	\item[(ii)] $\phi$ is  unimodal and symmetric about $1/2$, and
	\item[(iii)] $\phi(z) > 0$ if and only if $z \in (0,1)$,
\een
where $[x]$ denotes the largest integer less than or equal to $x$.
Figure \ref{fig:phi} shows an illustration of $\phi$ and related functions.
For a positive integer $m=m_n$, with $m_n \uparrow \infty$ as $n \to \infty$, let $z_j = j/m$, $I_j = [z_j, z_{j+1}]$ for $j=0, \ldots, m-1$, $J = \{0, 1, \ldots, m-1\}^d$ and $\phi_j(z) = \phi(m(z-z_j))$.
For a multi-index $\bj = (j_1, \ldots, j_d) \in J$ and $\alpha = (\alpha_\bj)_{\bj \in J} \in \{-1, +1\}^{|J|}$, define $\bg_\alpha: [0,1]^d \to \bbR^d$ as
\bean
	\bg_\alpha(\bz) = \left( z_1 + \frac{c_1}{m^\beta} \sum_{\bj \in J} \alpha_\bj \phi_{j_1}(z_1) \cdots \phi_{j_d}(z_d), z_2, \ldots, z_d \right),
\eean
where $c_1=c_1(\phi, d)$ is a small enough constant described below.
Then, it is easy to check that $\bg_\alpha$ is a one-to-one function from $[0,1]^d$ onto itself, and $\bg_\alpha \in \cH_K^\beta([0,1]^d) \times \cdots \times \cH_K^\beta([0,1]^d)$ for large enough $K=K(\beta, c_1)$.

Let  $\bZ = (Z_1, \ldots, Z_d)$ be a uniform random variable on $(0,1)^d$.
Then, by the change of variables formula, the Lebesgue density $q_\alpha$ of $\bY = \bg_\alpha(\bZ)$ is given as
\bean
	q_\alpha(\by) 
	= \left| \frac{\partial \bz}{\partial \by} \right|
	= \left( 1 + \frac{c_1}{m^{\beta}} \sum_{\bj \in J} \alpha_\bj \phi_{j_1}'(z_1) \phi_{j_2}(y_2) \cdots \phi_{j_d}(y_d) \right)^{-1}
\eean
for $\by \in [0,1]^d$, where $\phi'$ denotes the derivative of $\phi$.
Here, $z_1 = z_1(y_1, \ldots, y_d)$ is implicitly defined.

We first find an upper bound of $K(q_\alpha, q_{\alpha'})$ for $\alpha, \alpha' \in \{-1, +1\}^{|J|}$, where $K(p, q) = \int \log p / q dP$ is the Kullback--Leibler divergence.
Since $\beta \geq 1$, $q_\alpha$ is bounded from above and below for a small enough $c_1$.
Also, $\bg_\alpha(C_\bj) = C_\bj$, where $C_\bj = I_{j_1} \times \cdots \times I_{j_d}$.
Therefore, we have
\bean
	|q_\alpha (\by) - q_{\alpha'}(\by)|
	\lesssim \left|\frac{1}{q_\alpha (\by)} - \frac{1}{q_{\alpha'}(\by)} \right|
	\leq 2 \frac{c_1}{m^{\beta-1}} \|\phi'\|_\infty \|\phi\|_\infty^{d-1}.
\eean
Since the ratio $q_\alpha / q_{\alpha'}$ is bounded from above and below, we can use a well-known inequality $K(q_\alpha, q_{\alpha'}) \lesssim d_H^2 (q_\alpha, q_{\alpha'})$, where $d_H$ denotes the Hellinger distance; see Lemma B.2 from \cite{ghosal2017fundamentals}.
Since $|\sqrt{q_\alpha} - \sqrt{q_{\alpha'}}| \lesssim |q_\alpha - q_{\alpha'}|$, we have
\bean
	K(q_\alpha, q_{\alpha'})
	\lesssim \int_{[0,1]^d} | q_\alpha(\by) - q_{\alpha'}(\by) |^2 d \by
	\lesssim \frac{c_1^2 \|\phi'\|_\infty^2 \|\phi\|_\infty^{2(d-1)}}{m^{2(\beta-1)}}.
\eean

Next, we derive a lower bound for $W_1(q_\alpha, q_{\alpha'})$.
Suppose that $\alpha_\bj \neq \alpha_\bj'$ for some $\bj \in J$.
Then, the excess mass of $Q_\alpha$ over $Q_{\alpha'}$ on $C_\bj$ is
\bean
	&& \int_{\{\by\in C_j: q_\alpha(\by) > q_{\alpha'}(\by)\}} \Big\{q_\alpha (\by) - q_{\alpha'}(\by) \Big\} d \by
	\\
	&& =
	\frac{1}{2} \int_{C_\bj} |q_\alpha (\by) - q_{\alpha'}(\by)| d \by
	\\
	&& \gtrsim \int_{C_\bj}  \left|\frac{1}{q_\alpha (\by)} - \frac{1}{q_{\alpha'}(\by)} \right| d\by
	\\
	&& = \frac{2c_1}{m^\beta} \int_{C_\bj} | \phi_{j_1}'(z_1) \phi_{j_2}(y_2) \cdots \phi_{j_d}(y_d) | d\by
	\\
	&& = \frac{2c_1}{m^\beta} \int_{C_\bj} | \phi_{j_1}'(z_1) \phi_{j_2}(z_2) \cdots \phi_{j_d}(z_d) |  \left| \frac{\partial \by}{\partial \bz} \right| d\bz
	\\
	&& \gtrsim \frac{c_1}{m^{(\beta-1)+d}} \int_{(0,1)^d} |\phi'(z_1) \phi(z_2) \cdots \phi(z_d)| d\bz.
\eean
In virtue of Corollary 1.16 from \cite{villani2003topics}, with a (unique) optimal transport plan between $Q_\alpha$ and $Q_{\alpha'}$, some portion $\gamma \in (0,1)$ of this excess mass must be transported at least the distance of $c_2/m$, where constants $\gamma$ and $c_2$ can be chosen so that they depend only on $d$ and $\phi$.
Hence, for some constant $c_3 = c_3(\phi, d)$, 
\bean
	W_1(q_\alpha, q_{\alpha'}) \geq \frac{c_1 c_3}{m^{\beta+d}} H(\alpha, \alpha'),
\eean
where $H(\alpha, \alpha')= \sum_{\bj \in J} I(\alpha_\bj \neq \alpha'_\bj)$ denotes the Hamming distance between $\alpha$ and $\alpha'$.

With the Hamming distance on $\{-1, +1\}^{|J|}$, it is well-known (\eg\ see page 124 of \cite{wainwright2019high}) that there is a $|J|/4$-packing $\cA$ of $\{-1, +1\}^{|J|}$ whose cardinality is at least $e^{|J|/16}$.
Let $P_\alpha$ be the convolution of $Q_\alpha$ and $\cN(\bzero_d, \sigma_0^2 \Id_d)$.
Then, $K(p_\alpha, p_{\alpha'}) \leq K(q_\alpha, q_{\alpha'})$ by Lemma B.11 of \cite{ghosal2017fundamentals}.
By Fano's method (Proposition 15.12 from \cite{wainwright2019high}), we have
\bean
	\mathfrak{M}(\cG_0, \sigma_0)
	\gtrsim \frac{c_1 c_3}{m^\beta} \left\{ 1 - \frac{n c_1^2 C(\phi, d) m^{-2(\beta-1)} + \log 2 }{m^d/16} \right\}.
\eean
If $n \asymp m^{d + 2(\beta-1)}$, and $c_1$ is small enough,  we have the desired result.

\subsubsection*{Case 2: $q=0$ and $t_0 = d < D$}

Define a subset $\cG_1$ of $\cG_0 = \cH_K^\beta ([0,1]^d) \times \cdots \times \cH_K^\beta ([0,1]^d)$ as
\bean
    \cG_1 = \Big\{ \bg \in \cG_0: g_{d+1}(\bz) = \cdots = g_D(\bz) = 0 \Big\},
\eean
where $\bg(\cdot) = (g_1(\cdot), \ldots, g_D(\cdot))$.
The problem of obtaining a lower bound of the minimax risk $\mathfrak{M}(\cG_1, \sigma_0)$ reduces to Case 1, hence $\mathfrak{M}(\cG_1, \sigma_0)$ is bounded below by a multiple of $n^{-\beta / (2\beta+d-2)}$.
Since $\cG_1 \subset \cG_0$, we have $\mathfrak{M}(\cG_0, \sigma_0) \geq \mathfrak{M}(\cG_1, \sigma_0)$.

\subsubsection*{Case 3: $q=0$ and $t_0 < d \leq D$}

Similarly to Case 2, define a subset $\cG_2$ of $\cG_0$ as
\bean
    \cG_2 = \Big\{ \bg\in\cG_0 &:& \bg(\bz) = \big( g_1(\bz_{1:t_0}), \ldots, g_{t_0}(\bz_{1:t_0}), 0, \ldots, 0\big) 
    \\
    && \text{for some $g_j: [0,1]^{t_0} \to \bbR$, $j=1, \ldots, t_0$} \Big\},
\eean
where $\bz_{1:t_0} = (z_1, \ldots, z_{t_0})$.
Then, the problem reduces to Case 1 with $d$ replaced by $t_0$. Hence, we obtain a desired lower bound $\mathfrak{M}(\cG_0, \sigma_0) \geq \mathfrak{M}(\cG_2, \sigma_0) \gtrsim n^{-\beta / (2\beta+t_0-2)}$

\subsubsection*{Case 4: General $q$}

For $\cG_0 = \cG_0(q, \bd, \bt, \bbeta, K)$, fix $i_0 \in \{0, \ldots, q\}$.
We consider a subset $\cG_3$ of $\cG_0$ consisting of functions of the form $\bg = \bh_q \circ \bh_{q-1} \circ \cdots \circ \bh_1 \circ \bh_0$, where each $\bh_i: [a_i, b_i]^{d_i} \to [a_{i+1}, b_{i+1}]^{d_{i+1}}$ satisfies the following properties:
\ben
    \item For $i < i_0$, $\bh_i(\bx) = (\bx_{1:(d_i \wedge d_{i+1})}, \bzero_{d_{i+1} - d_i \wedge d_{i+1}})$.
    \item For $i > i_0$, $\bh_i(\bx) = (\bx_{1:t_{i_0}}, \bzero_{d_{i+1} - t_{i_0}})$.
    \item $\bh_{i_0}(\bx) = (h_{i_0 1}(\bx_{1:t_{i_0}}), \ldots, h_{i_0 t_{i_0}}(\bx_{1:t_{i_0}}), 0, \ldots, 0 )$ for some function $h_{i_0 j} \in \cH_K^{\beta_{i_0}}([a_{i_0}, b_{i_0}]^{t_{i_0}})$.
\een 
Since $t_{i_0} \leq \min\{d_0, \ldots, d_{q+1}\}$, we have
\bean
    \bg(\bz) = (h_{i_0 1}(\bz_{1:t_{i_0}}), \ldots, h_{i_0 t_{i_0}}(\bz_{1:t_{i_0}}), 0, \ldots, 0).
\eean
Again, the problem reduces to Case 1 with $(d, \beta)$ replaced by $(t_{i_0}, \beta_{i_0})$.
Therefore, $\mathfrak{M}(\cG_0, \sigma_0) \geq \mathfrak{M}(\cG_3, \sigma_0) \gtrsim n^{-\beta_{i_0} / (2\beta_{i_0} + t_{i_0} - 2)}$.
Since this inequality holds for all $i_0 \in \{0, \ldots, q\}$, the assertion of the theorem follows.
\qed

\end{appendix}

\end{document}